# Relative $K_0$, annihilators, Fitting ideals and the Stickelberger phenomena

Victor P. Snaith

**Abstract**

When $G$ is abelian and $l$ is a prime we show how elements of the relative K-group $K_0(\mathbf{Z}_l[G], \mathbf{Q}_l)$ give rise to annihilator/Fitting ideal relations of certain associated $\mathbf{Z}[G]$-modules. Examples of this phenomenon are ubiquitous. Particularly, we give examples in which $G$ is the Galois group of an extension of global fields and the resulting annihilator/Fitting ideal relation is closely connected to Stickelberger's Theorem and to the conjectures Coates-Sinnott and Brumer.

## 1 Introduction

**1.1** Suppose that $l$ is a prime and $G$ is a finite group. In this case the relative K-group $K_0(\mathbf{Z}_l[G], \mathbf{Q}_l)$ appearing in the localisation sequence (see §2.1) is isomorphic to the zero-th K-group of the category of finite $\mathbf{Z}_l[G]$-modules of finite projective dimension, usually denoted by $K_0T(\mathbf{Z}_l[G])$. When $G$ is abelian we have an isomorphism of the form $K_0(\mathbf{Z}_l[G], \mathbf{Q}_l) \cong \frac{\mathbf{Q}_l[G]^*}{\mathbf{Z}_l[G]^*}$.

As explained in §2, elements of $K_0(\mathbf{Z}_l[G], \mathbf{Q}_l)$ may be constructed from a bounded perfect complex of $\mathbf{Z}_l[G]$-modules together with a $\mathbf{Q}_l[G]$-trivialisation of its cohomology. In arithmetic and algebraic geometry there are many situations, some of which are studied in §§3-5, which give rise to this data. In this paper we shall consider the simplest example of this type of element, namely the case of a perfect complex of $\mathbf{Z}_l[G]$-modules whose cohomology groups are finite. Our main result concerns the simplest case of all, when the complex has only two non-zero, finite cohomology groups. In this case (Theorem 2.4) when $G$ is abelian, there are Stickelberger-type relations between the element of $\mathbf{Q}_l[G]^*$ which represents the relative $K_0$ class, the annihilators and the Fitting ideals of the cohomology groups. In order to apply Theorem 2.4 one must evaluate evaluate the element of $\mathbf{Q}_l[G]^*$ representing the class on $K_0(\mathbf{Z}_l[G], \mathbf{Q}_l)$. In general this is difficult to do but in the cases of abelian Galois extensions of function fields and totally real number fields this element turns out to be a classical higher order Stickelberger element (see Theorems 1.6 and 3.4).

Let us begin by recalling the higher Stickelberger elements and the conjectures pertaining to them. Suppose that $L/\mathbf{Q}$ is a finite Galois extension of number fields with abelian Galois group, $G(L/\mathbf{Q})$, and $L$ totally real. Then, for each integer $n \geq 2$, there is a unique unit of the rational group-ring

$$\Theta_{L/\mathbf{Q}}(n) \in \mathbf{Q}[G(L/\mathbf{Q})]^*$$

such that

$$\chi(\Theta_{L/\mathbf{Q}}(n)) = L_\mathbf{Q}(1-n, \chi^{-1})$$

for each one-dimensional complex representation $\chi$ where $L_\mathbf{Q}(s, \chi^{-1})$ denotes the Dirichlet L-function of the character $\chi^{-1}$ ([58] Ch.4). The rationality of $\Theta_{L/\mathbf{Q}}(n)$ is seen by writing the L-function in terms of partial zeta functions

$$L_\mathbf{Q}(1-n, \chi^{-1}) = \sum_{g \in G(L^{Ker(\chi)}/\mathbf{Q})} \chi(g)^{-1} \zeta_\mathbf{Q}(g, 1-n)$$

and recalling that $\zeta_\mathbf{Q}(g, 1-n)$ is a rational number, by a result of Klingen and Siegel (cf. [48]).

Inspired by Stickelberger's Theorem ([58] p.94), Coates and Sinnott made a fundamental conjecture concerning Galois actions on the algebraic K-groups of algebraic integers in number fields.

**Conjecture 1.2** ([11]; see also [12]).)
In the situation of §1.1, let $l$ be a prime and let $\mathcal{O}_L$ denote the ring of algebraic integers of $L$. For each $n \geq 1$

$$w_{n+1}(\mathbf{Q}) \Theta_{L/\mathbf{Q}}(n+1) \cdot \mathrm{ann}_{\mathbf{Z}_l[G(L/\mathbf{Q})]}(H^0_{\text{ét}}(\mathrm{Spec}(\mathcal{O}_L); \mathbf{Q}_l/\mathbf{Z}_l(n+1)))$$

$$\subseteq \mathrm{ann}_{\mathbf{Z}_l[G(L/\mathbf{Q})]}(K_{2n}(\mathcal{O}_L) \otimes \mathbf{Z}_l).$$

Here $w_{n+1}(\mathbf{Q})$ denotes the largest integer $m$ such that the Galois group, $G(\mathbf{Q}(e^{2\pi i/m})/\mathbf{Q})$, has exponent dividing $n+1$. Alternatively it is the order of $H^0_{\text{ét}}(\mathrm{Spec}(\mathbf{Z}); \mathbf{Q}/\mathbf{Z}(n+1))$ (see ([50] §7.2.4).

**Conjecture 1.3** In the situation of §1.1 and Conjecture 1.2 one could also conjecture the stronger result that, for each $n \geq 1$,

$$\Theta_{L/\mathbf{Q}}(n+1) \cdot \mathrm{ann}_{\mathbf{Z}_l[G(L/\mathbf{Q})]}(H^0_{\text{ét}}(\mathrm{Spec}(\mathcal{O}_L); \mathbf{Q}_l/\mathbf{Z}_l(n+1)))$$

$$\subseteq F_{\mathbf{Z}_l[G(L/\mathbf{Q})]}(K_{2n}(\mathcal{O}_L) \otimes \mathbf{Z}_l).$$

Here $F_{\mathbf{Z}_l[G(L/\mathbf{Q})]}(M)$ denotes the Fitting ideal of the module $M$ (see §2).

This conjecture is stronger than Conjecture 1.2 since the integer $w_{n+1}(\mathbf{Q})$ has been removed and the Fitting ideal is contained in the annihilator ideal (see §2.3).

**Remark 1.4** (i) The classical theorem of Stickelberger ([58] p.94) corresponds to the case when $n = 0$ in Conjecture 1.2 when $K_0(\mathcal{O}_L)$ is replaced by its torsion subgroup, the ideal class group. The Brumer Conjecture ([4], [60]), which is still open, corresponds to the case when $n = 0$ and $L/K$ is an abelian extension with $K$ totally real.

(ii) Furthermore, inspired by the Brumer Conjecture, a conjecture similar to Conjecture 1.3 makes sense for any abelian Galois extension of number fields, $L/K$, with $L$ totally real (see [50] Conjecture 7.2.6).

(iii) There is an equivariant Chern class homomorphism of the form (cf. [2] Theorem B; [14])

$$c_{n+1,2} : K_{2n}(\mathcal{O}_L) \otimes \mathbf{Z}_l \longrightarrow H^1_{\text{ét}}(\text{Spec}(\mathcal{O}_L[1/l]); \mathbf{Q}_l/\mathbf{Z}_l(n+1)).$$

The Lichtenbaum-Quillen Conjecture predicts (at least when $l$ is odd) that $c_{n+1,2}$ is an isomorphism. This was proved for $n = 1$ in [55]. As a corollary of the fundamental results of Voevodsky [56] [57], when $l = 2$ this Chern class is nearly an isomorphism in all dimensions [43]. Voevodsky's method requires the existence of suitable "norm varieties" which is not yet established for all odd primes. However, recent work by Rost combined with that of Suslin-Voevodsky shows that $c_{n+1,2}$ is an isomorphism for $l = 3$.

The annihilator

$$\text{ann}_{\mathbf{Z}_l[G(L/\mathbf{Q})]}(H^0_{\text{ét}}(\text{Spec}(\mathcal{O}_L); \mathbf{Q}_l/\mathbf{Z}_l(n+1)))$$

is well-known ([10] Lemma 2.3; [54] p.82). Explicitly, for $n \geq 0$, this annihilator is equal to the ideal generated by the elements, $(P, L/\mathbf{Q}) - NP^{n+1} \in \mathbf{Z}_l[G(L/\mathbf{Q})]$, where $P$ runs through primes which ramify in $L/\mathbf{Q}$ or divide

$$w_{n+1}(L) = |H^0_{\text{ét}}(\text{Spec}(\mathcal{O}_L); \mathbf{Q}_l/\mathbf{Z}_l(n+1))|.$$

Here $(P, L/\mathbf{Q}) \in G(L/\mathbf{Q})$ is the element corresponding to $P$ under Artin reciprocity.

In [14] it is shown that (the result can also be deduced from [15]) $K_{2n}(\mathcal{O}_L) \otimes \mathbf{Z}_l$ maps surjectively onto $H^1_{\text{ét}}(\text{Spec}(\mathcal{O}_L[1/l]); \mathbf{Q}_l/\mathbf{Z}_l(n+1))$ so that the Coates-Sinnott Conjecture of §1.2, as modified in Conjecture 1.3, predicts an inclusion of the form

$$\Theta_{L/\mathbf{Q}}(n+1) \cdot \text{ann}_{\mathbf{Z}_l[G(L/\mathbf{Q})]}(H^0_{\text{ét}}(\text{Spec}(\mathcal{O}_L); \mathbf{Q}_l/\mathbf{Z}_l(n+1)))$$

$$\subseteq F_{\mathbf{Z}_l[G(L/\mathbf{Q})]}(H^1_{\text{ét}}(\text{Spec}(\mathcal{O}_L[1/l]); \mathbf{Q}_l/\mathbf{Z}_l(n+1))).$$

(iv) The annihilator/Fitting ideal relation of (iii) makes sense more generally. Let $L/K$ be a Galois extension of totally real number fields with abelian Galois group, $G(L/K)$. Let $S$ be a finite set of primes of $K$ including all those which ramify in $L/K$. Then, as in §1.1, for each $n \geq 2$ there exists a unique unit of the rational group-ring

$$\Theta_{L/K,S}(n) \in \mathbf{Q}[G(L/K)]^*$$

such that

$$\chi(\Theta_{L/K,S}(n)) = L_{K,S}(1-n, \chi^{-1})$$

for each one-dimensional complex representation, $\chi$. Here $\chi$ is a character of $G(L/K)$ and $L_{K,S}(s, \chi^{-1})$ denotes the Artin L-function of $\chi^{-1}$ with the Euler factors associated to the primes in $S$ removed.

Let $S'$ be the set of places of $L$ over $S$. Set $X_l = \text{Spec}(O_{L,S'}[1/l])$, where $O_{L,S'}$ is the ring of $S'$-integers of $L$. For each $n \geq 1$, one would expect an annihilator relation of the form

$$\Theta_{L/K,S}(n+1) \cdot \text{ann}_{\mathbf{Z}_l[G(L/K)]}(H^0_{\text{ét}}(\text{Spec}(\mathcal{O}_L); \mathbf{Q}_l/\mathbf{Z}_l(n+1)))$$

$$\subseteq F_{\mathbf{Z}_l[G(L/K)]}(H^1_{\text{ét}}(X_l; \mathbf{Q}_l/\mathbf{Z}_l(n+1))) \triangleleft \mathbf{Z}_l[G(L/K)].$$

In fact, it can be shown that these ideals have the same radical ([50] Theorem 7.3.2).

(v) Using Iwasawa theory, it is shown in ([11] Theorem 2) that if $L/\mathbf{Q}$ is a totally real abelian extension of conductor $f$, $l$ is a prime and $b$ is an integer coprime to $fl$ then $w_2(\mathbf{Q})(b^2 - (b, L/\mathbf{Q}))\Theta_{L/\mathbf{Q}}(2) \in \mathbf{Q}_l[G(L/\mathbf{Q})]$ lies in

$$\text{ann}_{\mathbf{Z}_l[G(L/K)]}(H^1_{\text{ét}}(\text{Spec}(\mathcal{O}_L); \mathbf{Q}_l/\mathbf{Z}_l(2))) = \text{ann}_{\mathbf{Z}_l[G(L/K)]}(K_2(\mathcal{O}_L) \otimes \mathbf{Z}_l).$$

By combining ([11] Theorem 2.1) with results from [46], it is shown more generally in [39], under the same conditions (actually the condition that $(b, fl) = 1$ seems to have been omitted in ([39] Théorème 2.2)), that $w_{j+1}(\mathbf{Q})(b^{j+1} - (b, L/\mathbf{Q}))\Theta_{L/\mathbf{Q}}(j+1) \in \mathbf{Q}_l[G(L/\mathbf{Q})]$ lies in

$$\text{ann}_{\mathbf{Z}_l[G(L/K)]}(H^1_{\text{ét}}(\text{Spec}(\mathcal{O}_L); \mathbf{Q}_l/\mathbf{Z}_l(j+1)))$$

for $j = 1, 3, 5 \ldots$.

**1.5** Now let us describe the main results of this paper. Let $l$ be an odd prime and let $m$ be a positive integer prime to $l$. If $\xi_t$ denotes the root of unity $\xi_t = e^{2\pi i/t}$, consider the cyclotomic field $\mathbf{Q}(\xi_{ml^{s+1}})$ and its totally real subfield $\mathbf{Q}(\xi_{ml^{s+1}})^+$. Let $S_{ml}$ denote the primes of $\mathbf{Q}(\xi_{ml^{s+1}})^+$ which divide $ml$ and let $X_l^+ = \text{Spec}(\mathbf{Z}[\xi_{ml^{s+1}}]_{S_{ml}}^+)$ denote the spectrum of the $S_{ml}$-integers in $\mathbf{Q}(\xi_{ml^{s+1}})^+$. Note that $H^0_{\text{ét}}(X_l^+; \mathbf{Q}_l/\mathbf{Z}_l(1-r))$ is isomorphic to $H^0_{\text{ét}}(\mathbf{Z}[\xi_{ml^{s+1}}][1/ml]; \mathbf{Q}_l/\mathbf{Z}_l(1-r))$ for all negative odd integers $r$.

The following result will be proved in §4.3.

**Theorem 1.6**

Let $l$ be an odd prime. Let $l, m$ and $X_l^+$ be as in the notation of §1.5 and let $\Theta_{\mathbf{Q}(\xi_{ml^{s+1}})^+/\mathbf{Q}}(1-r)$ denote the higher Stickelberger element of §1.1. Then, for $r = -1, -3, -5, \ldots$ and any positive integer $s$:

(i) there exists a chain of annihilator ideal relations for étale cohomology of the form

$$\{t^{m_0} \mid t \in \mathrm{ann}_{\mathbf{Z}_l[G(\mathbf{Q}(\xi_{ml^{s+1}})^+/\mathbf{Q})]}(H^1_{\text{ét}}(X_l^+; \mathbf{Q}_l/\mathbf{Z}_l(1-r)))\}$$

$$\subseteq \Theta_{\mathbf{Q}(\xi_{ml^{s+1}})^+/\mathbf{Q}}(1-r)\mathrm{ann}_{\mathbf{Z}_l[G(\mathbf{Q}(\xi_{ml^{s+1}})^+/\mathbf{Q})]}(H^0_{\text{ét}}(X_l^+; \mathbf{Q}_l/\mathbf{Z}_l(1-r)))$$

$$\subseteq \mathrm{ann}_{\mathbf{Z}_l[G(\mathbf{Q}(\xi_{ml^{s+1}})^+/\mathbf{Q})]}(H^1_{\text{ét}}(X_l^+; \mathbf{Q}_l/\mathbf{Z}_l(1-r))).$$

(ii) If $l$ does not divide $m-1$ then in (i) the final

$$\mathrm{ann}_{\mathbf{Z}_l[G(\mathbf{Q}(\xi_{ml^{s+1}})^+/\mathbf{Q})]}(H^1_{\text{ét}}(X_l^+; \mathbf{Q}_l/\mathbf{Z}_l(1-r)))$$

may be replaced by

$$F_{\mathbf{Z}_l[G(\mathbf{Q}(\xi_{ml^{s+1}})^+/\mathbf{Q})]}(H^1_{\text{ét}}(X_l^+; \mathbf{Q}_l/\mathbf{Z}_l(1-r))).$$

Here $m_0$ is the minimal number of generators of the $\mathbf{Z}_l[G(\mathbf{Q}(\xi_{ml^{s+1}})^+/\mathbf{Q})]$-module $H^1_{\text{ét}}(X_l^+; \mathbf{Q}_l/\mathbf{Z}_l(1-r))$ and $F_{\mathbf{Z}_l[G(\mathbf{Q}(\xi_{ml^{s+1}})^+/\mathbf{Q})]}(M)$ denotes the Fitting ideal of $M$.

From the localisation sequence of ([52] p.268), when $r = -1, -3, -5, \ldots$, $H^1_{\text{ét}}(\mathrm{Spec}(\mathbf{Z}[\xi_{ml^{s+1}}]^+[1/l]); \mathbf{Q}_l/\mathbf{Z}_l(1-r))$ is a Galois submodule of $H^1_{\text{ét}}(X_l^+; \mathbf{Q}_l/\mathbf{Z}_l(1-r))$, which implies the following result.

**Corollary 1.7**

In the notation of §1.5 and Theorem 1.6, when $l$ is an odd prime not dividing $m$ and $r = -1, -3, -5, \ldots$

$$\Theta_{\mathbf{Q}(\xi_{ml^{s+1}})^+/\mathbf{Q}}(1-r)\mathrm{ann}_{\mathbf{Z}_l[G(\mathbf{Q}(\xi_{ml^{s+1}})^+/\mathbf{Q})]}(H^0_{\text{ét}}(\mathrm{Spec}(\mathbf{Z}[\xi_{ml^{s+1}}]^+[1/l]); \mathbf{Q}_l/\mathbf{Z}_l(1-r)))$$

$$\subseteq \mathrm{ann}_{\mathbf{Z}_l[G(\mathbf{Q}(\xi_{ml^{s+1}})^+/\mathbf{Q})]}(H^1_{\text{ét}}(\mathrm{Spec}(\mathbf{Z}[\xi_{ml^{s+1}}]^+[1/l]; \mathbf{Q}_l/\mathbf{Z}_l(1-r))))$$

for all positive integers $s$.

The paper is organised in the following manner. In §2 we introduce the relative K-group $K_0(\mathbf{Z}_l[G], \mathbf{Q}_l)$ for any finite group $G$ and explain how elements of this group can yield annihilator/Fitting ideal relations (Theorem 2.4). Two important classes of elements to which the construction of §2 applies arise in connection with the étale cohomology groups of curves over finite fields and rings of integers in number fields. In §3 we apply Theorem 2.4 to an example from ([8] §7) to prove a Coates-Sinnott type of result (Theorem 3.4) for the étale cohomology groups of curves over finite fields. In §4 we use an example from [5] to prove a similar result (Theorem 1.6 proved in §4.3) for the étale cohomology $H^1_{\text{ét}}(X_l^+; \mathbf{Q}_l/\mathbf{Z}_l(1-r))$ of the S-integers in the totally real subfield of a cyclotomic field. When $r = -2, -4, -6, \ldots$ – the case not covered by Theorem 1.6 – we apply Theorem 2.4 in Theorem 4.6 to study the annihilator ideal of $H^2_{\text{ét}}(X_l^+; \mathbf{Z}_l(1-r))$. In §4.8 and Theorem 4.9 we use Theorem 4.6 together with calculation of Beilinson [3] to construct elements of this annihilator from the leading terms of the Dirichlet L-function at $s = r$. In §5 we discuss invariants lying in $K_0(\mathbf{Z}_l[G], \mathbf{Q}_l)$ which are constructed from relatively abelian extensions of totally real number fields with vanishing Iwasawa $\mu$-invariants, from the Galois action on vanishing cycles and from étale coverings of curves and surfaces.

Finally, a remark about the classical Stickelberger Theorem and Brumer Conjecture which, as mentioned in Remark 1.4, correspond to the case when $r = 0$. In ([50] Remark 7.2.11(ii)) it is explained how one might use "Tate sequence" constructed in [41] to approach the Brumer Conjecture by the method used to prove Theorem 1.6. This approach is much more difficult than the proof of Theorem 1.6 because of the presence of non-trivial regulators (cf. [50] §7.1.12). On the other hand, this approach would aim to prove a statement (analogous to that of Theorem 1.6) stating that at odd primes $l$, the Stickelberger ideal lies in the annihilator (and sometimes the Fitting ) ideal of the class-group of the $S$-integers (where $S$-contains all ramified primes). Even in the case of the classical Stickelberger Theorem this seems to lead to questions which are at present unanswered. In the classical Stickelberger case such a Fitting ideal inclusion fails for $l = 2$ but holds at odd primes when, for example, both the roots of unity and the class-group are cohomologically trivial.

Very recently I became aware of the very important results of Masato Kurihara [27] which show that the Fitting ideal appearing in Theorem 1.6(ii) is given by a generalised Stickelberger ideal. Using Kurihara's results one may verify the second inclusion in the statements of Theorems 1.6(ii) *without any restrictions on the integer m* while the other inclusion may in turn be combined with the results of [27] to yield upper bounds on the number of $\mathbf{Z}_l[G]$-module generators of the étale cohomology groups.

**Ackowledgments** The results of [5] are crucial to this paper and I am particularly indebted to David Burns for bringing his work to my attention. This work was begun during a visit to the University of Pennsylvania in March 2000 and completed during a visit, funded by a Special Project Grant from the Royal Society, to the Euler Institute in St Petersburg during September 2000. I am particularly grateful to these institutions for their hospitality and to mathematicians there – Ted Chinburg, Ivan Panin and Igor Zhukov

– for help and encouragement. In addition, conversations with John Coates, Cornelius Greither, Bernhard Köck, Masato Kurihara, Christian Popescu, Nguyen Thong Quang Do, Jürgen Ritter, Al Weiss and Andrew Wiles have been very helpful. In fact, [11] was one of the papers which first got me interested in algebraic K-theory a long time ago.

## 2 $K_0(\mathbf{Z}_l[G], \mathbf{Q}_l)$ and annihilator/Fitting ideal relations

**2.1** Let $l$ be a prime, $G$ a finite group and let $f : \mathbf{Z}_l[G] \longrightarrow \mathbf{Q}_l[G]$ denote the homomorphism of group-rings induced by the inclusion of the $l$-adic integers into the fraction field, the $l$-adic rationals. Write $K_0(\mathbf{Z}_l[G], \mathbf{Q}_l)$ for the relative K-group of $f$, denoted by $K_0(\mathbf{Z}_l[G], f)$ in ([53] p.214; see also [50] Definition 2.1.5). By ([53] Lemma 15.6) elements of $K_0(\mathbf{Z}_l[G], \mathbf{Q}_l)$ are represented by triples $[A, g, B]$ where $A, B$ are finitely generated, projective $\mathbf{Z}_l[G]$-modules and $g$ is a $\mathbf{Q}_l[G]$-module isomorphism of the form $g : A \otimes_{\mathbf{Z}_l} \mathbf{Q}_l \stackrel{\cong}{\to} B \otimes_{\mathbf{Z}_l} \mathbf{Q}_l$. Defining an exact sequence of triples in the obvious manner, the relations between these elements are generated by the following two types:

(i)  $[A, g, B] = [A', g', B'] + [A'', g'', B'']$ if there exists an exact sequence

$$0 \longrightarrow (A', g', B') \longrightarrow (A, g, B) \longrightarrow (A'', g'', B'') \longrightarrow 0$$

and

(ii)  $[A, gh, C] = [A, h, B] + [B, g, C]$.

This group fits into a localisation sequence of the form ([44] §5 Theorem 5; [20] p.233)

$$K_1(\mathbf{Z}_l[G]) \stackrel{f_*}{\longrightarrow} K_1(\mathbf{Q}_l[G]) \stackrel{\partial}{\longrightarrow} K_0(\mathbf{Z}_l[G], \mathbf{Q}_l) \stackrel{\pi}{\longrightarrow} K_0(\mathbf{Z}_l[G]) \stackrel{f_*}{\longrightarrow} K_0(\mathbf{Q}_l[G]).$$

Assume now that $G$ is abelian. In this case $K_1(\mathbf{Q}_l[G]) \cong \mathbf{Q}_l[G]^*$ because $\mathbf{Q}_l[G]$ is a product of fields and $K_1(\mathbf{Z}_l[G]) \cong \mathbf{Z}_l[G]^*$ ([13]I p.179 Theorem (46.24)). Under these isomorphisms $f_*$ is identified with the canonical inclusion.

The homomorphism, $K_0(\mathbf{Z}_l[G]) \stackrel{f_*}{\longrightarrow} K_0(\mathbf{Q}_l[G])$, is injective for all finite groups $G$ ([47] Theorem 34 p.131; [13]II p.47 Theorem 39.10). Alternatively, when $G$ is abelian, $\mathbf{Z}_l[G]$ is semi-local and the injectivity of $f_*$ follows from the fact that a finitely generated projective module over a local ring is free ([16] p.124 Corollary 4.8 and p.205 Exercise 7.2). Thus the localisation sequence yields an isomorphism of the form

$$K_0(\mathbf{Z}_l[G], \mathbf{Q}_l) \cong \frac{\mathbf{Q}_l[G]^*}{\mathbf{Z}_l[G]^*}$$

when $G$ is abelian. From the explicit description of $\partial$ ([53] p.216) this isomorphism sends the coset of $\alpha \in \mathbf{Q}_l[G]^*$ to $[\mathbf{Z}_l[G], (\alpha \cdot -), \mathbf{Z}_l[G]]$. The inverse isomorphism sends $[A, g, B]$, where $A$ and $B$ may be assumed to be free $\mathbf{Z}_l[G]$-modules, to the coset of $det(g) \in \mathbf{Q}_l[G]^*$ with respect to any choice of $\mathbf{Z}_l[G]$-bases for $A$ and $B$.

**Example 2.2** We shall be particularly interested in the following source of elements of $K_0(\mathbf{Z}_l[G], \mathbf{Q}_l)$.

As in §2.1, let $l$ be a prime and let $G$ be a finite abelian group. Suppose that
$$0 \longrightarrow F_k \xrightarrow{d_k} F_{k-1} \xrightarrow{d_{k-1}} \ldots \xrightarrow{d_2} F_1 \xrightarrow{d_1} F_0 \longrightarrow 0$$
is a bounded complex of finitely generated, projective $\mathbf{Z}_l[G]$-modules (i.e. a *perfect* complex of $\mathbf{Z}_\mathbf{l}[G]$-modules), having all its homology groups finite.

As usual, let $Z_t = \mathrm{Ker}(d_t : F_t \longrightarrow F_{t-1})$ and $B_t = d_{t+1}(F_{t+1}) \subseteq F_t$ denote the $\mathbf{Z}_l[G]$-modules of $t$-dimensional cycles and boundaries, respectively. We have short exact sequences of the form
$$0 \longrightarrow B_i \xrightarrow{\phi_i} Z_i \longrightarrow H_i(F_*) \longrightarrow 0$$
and
$$0 \longrightarrow Z_{i+1} \xrightarrow{\psi_{i+1}} F_{i+1} \xrightarrow{d_{i+1}} B_i \longrightarrow 0.$$

Applying $(- \otimes \mathbf{Q}_l)$ we obtain isomorphisms
$$\phi_i : B_i \otimes \mathbf{Q}_l \xrightarrow{\cong} Z_i \otimes \mathbf{Q}_l$$
and we may choose $\mathbf{Q}_l[G]$-module splittings of the form
$$\eta_i : B_i \otimes \mathbf{Q}_l \longrightarrow F_{i+1} \otimes \mathbf{Q}_l$$
such that $(d_{i+1} \otimes 1)\eta_i = 1 : B_i \otimes \mathbf{Q}_l \longrightarrow B_i \otimes \mathbf{Q}_l$.

Then we form a $\mathbf{Q}_l[G]$-module isomorphism of the form
$$X : \oplus_j F_{2j} \otimes \mathbf{Q}_l \xrightarrow{\cong} \oplus_j F_{2j+1} \otimes \mathbf{Q}_l$$
given by the composition
$$\oplus_j F_{2j} \otimes \mathbf{Q}_l \xrightarrow{\oplus_j (\psi_{2j} + \eta_{2j-1})^{-1}} \oplus_j (Z_{2j} \otimes \mathbf{Q}_l) \oplus (B_{2j-1} \otimes \mathbf{Q}_l)$$
$$\xrightarrow{\oplus_j (\phi_{2j}^{-1} \oplus \phi_{2j-1})} \oplus_j (B_{2j} \otimes \mathbf{Q}_l) \oplus (Z_{2j-1} \otimes \mathbf{Q}_l)$$
$$\xrightarrow{(\oplus_j \ \eta_{2j} + \psi_{2j-1})} \oplus_j F_{2j-1} \otimes \mathbf{Q}_l.$$

This construction defines a class, $[\oplus_j F_{2j}, X, \oplus_j F_{2j+1}]$, in $K_0(\mathbf{Z}_l[G], \mathbf{Q}_l)$ which is well-known to be independent of the choices of the splittings used to define $X$ ([53] Ch. 15; see also [50] Propositions 2.5.35 and 7.1.8).

We shall denote by
$$det(X) \in \frac{\mathbf{Q}_l[G]^*}{\mathbf{Z}_l[G]^*}$$
the element which corresponds to $[\oplus_j F_{2j}, X, \oplus_j F_{2j+1}] \in K_0(\mathbf{Z}_l[G], \mathbf{Q}_l)$ under the isomorphism of §2.1.

**2.3** Recall that, if $R$ is a ring and $M$ a (left) $R$-module, the (left) annihilator ideal $\text{ann}_R(M) \triangleleft R$ of $M$ is defined to be
$$\text{ann}_R(M) = \{r \in R \mid r \cdot m = 0 \text{ for all } m \in M\}.$$

Let us recall from ([35] Appendix; see also [59]) the properties of the Fitting ideal (referred to as the initial Fitting invariant in [40]).

Let $R$ be a commutative ring with identity and let $M$ be a finitely presented $R$-module, in our applications $M$ will actually be finite. Suppose that $M$ has a presentation of the form
$$R^a \xrightarrow{f} R^b \longrightarrow M \longrightarrow 0$$
with $a \geq b$ then the Fitting ideal of the $R$-module $M$, denoted by $F_R(M)$, is the ideal of $R$ generated by all $b \times b$ minors of any matrix representing $f$.

The Fitting ideal $F_R(M)$ is independent of the presentation chosen for $M$ and is contained in the annihilator ideal of $M$, $F_R(M) \subseteq \text{ann}_R(M)$. If $M$ is generated by $n$ elements then $\text{ann}_R(M)^n \subseteq F_R(M)$ and if $\pi : M \to M'$ is a surjection of finitely presented $R$-modules then $F_R(M) \subseteq F_R(M')$.

The following result yields relations between the annihilator ideals and Fitting ideals of the homology modules in Example 2.2 in the special case when each $H_i(F_*)$ is finite and zero except for $i = 0, 1$.

**Theorem 2.4** *Let $G$ be a finite abelian group and $l$ a prime. Suppose that*
$$0 \longrightarrow F_k \xrightarrow{d_k} F_{k-1} \xrightarrow{d_{k-1}} \ldots \xrightarrow{d_2} F_1 \xrightarrow{d_1} F_0 \longrightarrow 0$$
*is a bounded, perfect complex of $\mathbf{Z}_l[G]$-modules, as in Example 2.2, having $H_i(F_*)$ finite for $i = 0, 1$ and zero otherwise. Let*
$$[\oplus_j F_{2j}, X, \oplus_j F_{2j+1}] \in K_0(\mathbf{Z}_l[G], \mathbf{Q}_l) \cong \frac{\mathbf{Q}_l[G]^*}{\mathbf{Z}_l[G]^*}$$
*be as in Example 2.2. Then:*

(i) if $t_i \in \mathrm{ann}_{\mathbf{Z}_l[G]}(H_i(F_*))$,

$$det(X)^{(-1)^i} t_i^{m_i} \in \mathrm{ann}_{\mathbf{Z}_l[G]}(H_{1-i}(F_*)) \lhd \mathbf{Z}_l[G]$$

for $i = 0, 1$. Here $m_0, m_1$ is the minimal number of generators required for the $\mathbf{Z}_l[G]$-module $H_0(F_*), Hom(H_1(F_*), \mathbf{Q}_l/\mathbf{Z}_l)$, respectively,

(ii) if the Sylow $l$-subgroup of $G$ is cyclic then in (i) $\mathrm{ann}_{\mathbf{Z}_l[G]}(H_{1-i}(F_*))$ may be replaced by $F_{\mathbf{Z}_l[G]}(H_{1-i}(F_*))$.

**2.5** The proof of Theorem 2.4 will occupy the remainder of this section, culminating in §2.11. We begin with some preliminary observations and constructions.

In the situation of Theorem 2.4, without changing the class $[\oplus_j F_{2j}, X, \oplus_j F_{2j+1}]$, we may replace $F_0, \ldots, F_{k-1}$ by finitely generated, free $\mathbf{Z}_l[G]$-modules. Then, since $\sum_i (-1)^i [F_i] = f_*(\pi([\oplus_j F_{2j}, X, \oplus_j F_{2j+1}])) = 0$ and $f_*$ is injective, $F_k$ must be stably free and hence free. Therefore henceforth we shall assume that each $F_i$ is a finitely generated, free $\mathbf{Z}_l[G]$-module. In this case we are going to modify the isomorphism $X$ to make it $\mathbf{Z}_l$-integral.

Firstly we observe that, for $i \geq 1$, the short exact sequence of $\mathbf{Z}_l[G]$-module homomorphisms

$$0 \longrightarrow B_{i+1} \xrightarrow{\psi_{i+1}} F_{i+1} \xrightarrow{d_{i+1}} B_i \longrightarrow 0$$

splits. To see this, write $G = H \times G_1$ where $G_1$ is the Sylow $l$-subgroup of $G$ and observe, by induction starting with $B_{k-1} \cong F_k$, that $B_i$ is a torsion-free, cohomologically trivial $\mathbf{Z}_l[G_1]$-module if $i \geq 1$. Hence, by ([1] Theorem 7) there exist elements $v_1, \ldots, v_t \in B_i$ such that the natural map induces an isomorphism of the form

$$\mathbf{Z}/l^m[G_1] < v_1, \ldots, v_t > \xrightarrow{\cong} B_i/l^m B_i$$

for each $m \geq 1$. Since $B_i$ is $l$-adically complete, the inverse limit of these isomorphisms shows that $B_i$ is a free $\mathbf{Z}_l[G_1]$-module. Hence we may split $d_{i+1}$ as a $\mathbf{Z}_l[G_1]$-module homomorphism and, by averaging over $H$, as a $\mathbf{Z}_l[G]$-module homomorphism also.

Therefore we have

$$\eta_i : B_i \longrightarrow F_{i+1}$$

for $i \geq 1$ and

$$\phi_i : B_i \xrightarrow{=} Z_i$$

for $i \geq 2$ without tensoring with the $l$-adic rationals.

Since $H_0(F_*)$ is finite, in the notation of Example 2.2, we have a homomorphism

$$\eta = \eta_0 \cdot (\phi_0 \otimes 1)^{-1} : F_0 \otimes \mathbf{Q}_l \longrightarrow F_1 \otimes \mathbf{Q}_l$$

such that $(d_1 \otimes 1) \cdot \eta = 1$ on $F_0 \otimes \mathbf{Q}_l$. By definition, $X(F_0 \otimes \mathbf{Q}_l) \subseteq F_1 \otimes \mathbf{Q}_l \subseteq \oplus_j F_{2j-1} \otimes \mathbf{Q}_l$ and $X$ restricted to $F_0 \otimes \mathbf{Q}_l$ is equal to $\eta$.

Let $\mathbf{Q}_l[G]^*$ denote the units of $\mathbf{Q}_l[G]$ and let $t \in \mathbf{Z}_l[G] \cap \mathbf{Q}_l[G]^*$ satisfy $tH_0(F_*) = 0$. We may construct $\eta$ to satisfy $\eta(tF_0) \subseteq F_1$. For each free generator, $z \in F_0$, we have $tz = d_1(w)$ and so we may define $\eta(z) = t^{-1}w$. Since $G$ is abelian, multiplication by $t$ gives a $\mathbf{Z}_l[G]$-module homomorphism

$$\tilde{T} = (t \cdot -) : F_0 \longrightarrow F_0$$

and we may form the composition

$$X \cdot (\tilde{T} \oplus 1 \oplus 1 \oplus \ldots) : \oplus_j F_{2j} \otimes \mathbf{Q}_l \xrightarrow{\cong} \oplus_j F_{2j+1} \otimes \mathbf{Q}_l.$$

Our next result will show that $X \cdot (\tilde{T} \oplus 1 \oplus 1 \oplus \ldots)$ maps $\oplus_j F_{2j}$ to $\oplus_j F_{2j+1}$. To prove this it will be useful to have the explicit formula for $X$ on $(w_0, w_2, \ldots) \in \oplus_j F_{2j}$, namely

$$X(w_0, w_2, \ldots)$$

$$= (\eta(w_0) + d_2(w_2), \eta_2(w_2 - \eta_1(d_2(w_2))) + d_4(w_4), \ldots$$

$$\eta_{2t}(w_{2t} - \eta_{2t-1}(d_{2t}(w_{2t}))) + d_{2t+2}(w_{2t+2}), \ldots).$$

**Lemma 2.6**

Let $X$, $\tilde{T} = (t \cdot -)$ and $t \in \mathbf{Q}_l[G]^*$ be as in §2.5. Then the $\mathbf{Q}_l[G]$-module isomorphism

$$X \cdot (\tilde{T} \oplus 1 \oplus 1 \oplus \ldots) : \oplus_j F_{2j} \otimes \mathbf{Q}_l \xrightarrow{\cong} \oplus_j F_{2j+1} \otimes \mathbf{Q}_l$$

restricts to yield an injective $\mathbf{Z}_l[G]$-module homomorphism of the form

$$X \cdot (\tilde{T} \oplus 1 \oplus \ldots \oplus 1) : \oplus_j F_{2j} \longrightarrow \oplus_j F_{2j-1}.$$

**Proof**

The injectivity follows from the fact that $X$ and $\tilde{T}$ are both $\mathbf{Q}_l[G]$-module isomorphisms. For the rest it suffices to verify that $X \cdot (\tilde{T} \oplus 1 \oplus \ldots \oplus 1)(F_{2s}) \subseteq \oplus_j F_{2j-1}$ for each $s \geq 0$. If $s = 0$ this is assured by the previous observation that $X$ equals $\eta$ on $F_0$. For $s \geq 1$, the integrality follows from the explicit formula of §2.5. □

**Proposition 2.7**

Suppose that $F_*$ is a bounded complex of free $\mathbf{Z}_l[G]$-modules with $X$, $\tilde{T} = (t \cdot -)$ and $t \in \mathbf{Q}_l[G]^* \cap \mathrm{ann}_{\mathbf{Z}_l[G]}(H_0(F_*))$ as in §2.5. Then

$$det(X) t^{rank_{\mathbf{Z}_l[G]}(F_0)} \in \mathrm{ann}_{\mathbf{Z}_l[G]}(H_1(F_*)) \triangleleft \mathbf{Z}_l[G]$$

and, if the Sylow $l$-subgroup of $G$ is cyclic, the annihilator $\mathrm{ann}_{\mathbf{Z}_l[G]}(H_1(F_*))$ may be replaced by the Fitting ideal $F_{\mathbf{Z}_l[G]}(H_1(F_*))$.

**Proof**

The short exact sequence of $\mathbf{Z}_l[G]$-modules

$$0 \longrightarrow \oplus_j F_{2j} \xrightarrow{X \cdot (\tilde{T} \oplus 1 \oplus \ldots)} \oplus_j F_{2j+1} \longrightarrow \frac{\oplus_j F_{2j+1}}{X \cdot (\tilde{T} \oplus 1 \oplus \ldots)(\oplus_j F_{2j})} \longrightarrow 0$$

implies, by definition of the Fitting ideal, that

$$det(X) t^{rank_{\mathbf{Z}_l[G]}(F_0)} = (det(X \cdot (\tilde{T} \oplus 1 \oplus \ldots)) \in F_{\mathbf{Z}_l[G]}(\frac{\oplus_j F_{2j+1}}{X \cdot (\tilde{T} \oplus 1 \oplus \ldots)(\oplus_j F_{2j})}).$$

Recall that the Fitting ideal is a submodule of the annihilator ideal. Now consider the $\mathbf{Z}_l[G]$-module

$$N_1 = \frac{\text{Ker}(d_1)}{\text{Ker}(d_1) \cap (X \cdot (\tilde{T} \oplus 1 \oplus \ldots))(\oplus_j F_{2j})} \cong \frac{\text{Ker}(d_1)(X \cdot (\tilde{T} \oplus 1 \oplus \ldots))(\oplus_j F_{2j})}{(X \cdot (\tilde{T} \oplus 1 \oplus \ldots))(\oplus_j F_{2j})}$$

which is a finite group, since it is a submodule of $N_2 = \frac{\oplus_j F_{2j+1}}{X \cdot (\tilde{T} \oplus 1 \oplus \ldots)(\oplus_j F_{2j})}$.

For $i = 1, 2$ define two finite $\mathbf{Z}_l[G]$-modules by $M_i = Hom(N_i, \mathbf{Q}_l/\mathbf{Z}_l)$ then there is a surjection of the form $\pi : M_2 \longrightarrow M_1$ which is the Pontrjagin dual of the inclusion homomorphism. Hence $\text{ann}_{\mathbf{Z}_l[G]}(M_2) \subseteq \text{ann}_{\mathbf{Z}_l[G]}(M_1)$. By ([35] Appendix 1) there is a similar inclusion for Fitting ideals $F_{\mathbf{Z}_l[G]}(M_2) \subseteq F_{\mathbf{Z}_l[G]}(M_1)$. Let $\tau$ denote the involution of the group-ring induced by sending each element of $G$, which is abelian, to its inverse. However, if $N$ is finite, the action of $G$ on $Hom(N, \mathbf{Q}_l/\mathbf{Z}_l)$ is given by $g(f) = f(g^{-1} \cdot -)$ which may be identified with $G$ acting on $N$ with the new action, $g(z) = \tau(g) \cdot z$. Therefore we may apply $\tau$ to obtain

$$\text{ann}_{\mathbf{Z}_l[G]}(N_2) = \tau(\text{ann}_{\mathbf{Z}_l[G]}(M_2)) \subseteq \tau(\text{anna}_{\mathbf{Z}_l[G]}(M_1)) = \text{ann}_{\mathbf{Z}_l[G]}(N_1)$$

in part (i). By Proposition 2.8, since in part (ii) the Sylow $l$-subgroup of $G$ is assumed to be cyclic, there is a similar chain of relations for Fitting ideals

$$F_{\mathbf{Z}_l[G]}(N_2) = \tau(F_{\mathbf{Z}_l[G]}(M_2)) \subseteq \tau(F_{\mathbf{Z}_l[G]}(M_1)) = F_{\mathbf{Z}_l[G]}(N_1).$$

On the other hand, let $w_j \in F_j$ and suppose that $(w_1, 0, 0, \ldots) = X \cdot (\tilde{T} \oplus 1 \oplus \ldots)(w_0, w_2, \ldots)$ with $w_1 \in Ker(d_1) \subseteq F_1$. The formula of §2.5 implies that $w_1 = \eta(\tilde{T}(w_0)) + d_2(w_2)$. Applying $d_1$ we obtain $0 = d_1(\eta(\tilde{T}(w_0)) + d_2(w_2)) = \tilde{T}(w_0)$, which yields $w_1 = d_2(w_2) \in d_2(F_2)$. Therefore $H_1(F_*) = \text{Ker}(d_1)/d_2(F_2)$ is a quotient of $N_1$ so that we have inclusions $\text{ann}_{\mathbf{Z}_l[G]}(N_1) \subseteq \text{ann}_{\mathbf{Z}_l[G]}(H_1(F_*))$ and $F_{\mathbf{Z}_l[G]}(N_1) \subseteq F_{\mathbf{Z}_l[G]}(H_1(F_*))$, by ([35] Appendix). This discussion yields

$$det(X) t^{rank_{\mathbf{Z}_l[G]}(F_0)} \in \text{ann}_{\mathbf{Z}_l[G]}(N_2) \subseteq \text{ann}_{\mathbf{Z}_l[G]}(N_1) \subseteq \text{ann}_{\mathbf{Z}_l[G]}(H_1(F_*)),$$

in part (i) and

$$det(X) t^{rank_{\mathbf{Z}_l[G]}(F_0)} \in F_{\mathbf{Z}_l[G]}(N_2) \subseteq F_{\mathbf{Z}_l[G]}(N_1) \subseteq F_{\mathbf{Z}_l[G]}(H_1(F_*)),$$

in part (ii), as required. $\square$

**Proposition 2.8** *(Dualisation and Fitting ideals)*

Let $l$ be a prime, $G$ a finite abelian group and $M$ a finite $\mathbf{Z}_l[G]$-module. Let $\tau$ be the involution sending each $g \in G$ to its inverse. Then

$$F_{\mathbf{Z}_l[G]}(M) = \tau(F_{\mathbf{Z}_l[G]}(Hom(M, \mathbf{Q}_l/\mathbf{Z}_l)))$$

for all finite $M$ if and only if the Sylow $l$-subgroup is cyclic.

**Proof**

The following argument, which was shown to me by Al Weiss, shows that whenever $G$ is an abelian group whose $l$-Sylow subgroup is non-cyclic then there exists a finite $\mathbf{Z}_l[G]$-module $M$ for which

$$F_{\mathbf{Z}_l[G]}(M) \neq \tau(F_{\mathbf{Z}_l[G]}(Hom(M, \mathbf{Q}_l/\mathbf{Z}_l)))$$

where $\tau$ is the involution of the group-ring induced by sending a group element to its inverse. Here the action on $Hom(M, \mathbf{Q}_l/\mathbf{Z}_l)$ is given by $g(f) = f(g^{-1} \cdot -)$ as in the proof of Proposition 2.7.

If the $l$-Sylow subgroup of $G$ is non-cyclic choose a surjection of the form

$$\pi : \mathbf{Z}_l[G] \longrightarrow \mathbf{Z}_l[C_l \times C_l]$$

where $C_l$ denotes the cyclic group of order $l$. Set $I = Ker(\pi)$. By the ([35] Appendix 4) we have

$$F_{\mathbf{Z}_l[C_l \times C_l]}(M/I \cdot M) = F_{\mathbf{Z}_l[G]/I}(M/I \cdot M) = Im(F_{\mathbf{Z}_l[G]}(M) \xrightarrow{\pi} \mathbf{Z}_l[C_l \times C_l]).$$

Now set $M = Ker(\mathbf{F}_l[C_l \times C_l] \xrightarrow{\epsilon} \mathbf{F}_l)$, the kernel of the modulo $l$ augmentation homomorphism. Hence

$$F_{\mathbf{Z}_l[C_l \times C_l]}(M/I \cdot M) = F_{\mathbf{Z}_l[C_l \times C_l]}(M) = Im(F_{\mathbf{Z}_l[G]}(M) \xrightarrow{\pi} \mathbf{Z}_l[C_l \times C_l])$$

and

$$F_{\mathbf{Z}_l[C_l \times C_l]}(M^\#/I \cdot M^\#) = F_{\mathbf{Z}_l[C_l \times C_l]}(M^\#) = Im(F_{\mathbf{Z}_l[G]}(M^\#) \xrightarrow{\pi} \mathbf{Z}_l[C_l \times C_l])$$

where $M^\# = Hom(M, \mathbf{Q}_l/\mathbf{Z}_l)$. We shall show that these two ideals are different.

Set $J = Ker(\mathbf{Z}_l[C_l \times C_l] \xrightarrow{\epsilon} \mathbf{Z}_l) \lhd \mathbf{Z}_l[C_l \times C_l]$, the augmentation ideal.

Firstly $M/J \cdot M$ is equal to the $C_l \times C_l$-coinvariants of $M$ and this is isomorphic to $C_l \times C_l$ as an abelian group. Hence, by ([35] Appendix 4 and 5),

$$l^2 \mathbf{Z}_l = F_{\mathbf{Z}_l}(C_l \times C_l) = F_{\mathbf{Z}_l[C_l \times C_l]/J}(M/J \cdot M) = Im(\epsilon : F_{\mathbf{Z}_l[C_l \times C_l]}(M) \longrightarrow \mathbf{Z}_l).$$

Secondly the $C_l \times C_l$-coinvariants of $M^\#$ are equal to $(M^{C_l \times C_l})^\# \cong \mathbf{F}_l$ with the trivial action so that

$$l\mathbf{Z}_l = F_{\mathbf{Z}_l}(\mathbf{F}_l) = F_{\mathbf{Z}_l[C_l \times C_l]/J}(M^\#/J \cdot M^\#) = Im(\epsilon : F_{\mathbf{Z}_l[C_l \times C_l]}(M^\#) \longrightarrow \mathbf{Z}_l),$$

which completes the first half of the proof.

Now suppose that the Sylow $l$-subgroup of $G$ is cyclic. The $l$-adic group-ring $\mathbf{Z}_l[G]$ is finite product of group-rings of the form $R[G_1]$ where $G_1$ is the the Sylow $l$-subgroup of $G$ and $R$ is the ring of integers in an unramified extension of $\mathbf{Q}_l$. In fact, $\mathbf{Z}_l[G] \cong \prod_\chi \mathbf{Z}_l[\chi][G_1]$ where $\chi$ runs through the Galois orbits of those characters of $G$ whose order is prime to $l$ ([49] §7.3.6; [59], [60]). Therefore $R[G_1]$ is a quotient of the polynomial ring $R[X]$ and the relation

$$F_{\mathbf{Z}_l[G]}(M) = \tau(F_{\mathbf{Z}_l[G]}(Hom(M, \mathbf{Q}_l/\mathbf{Z}_l)))$$

for all finite $\mathbf{Z}_l[G]$-module $M$ follows from the corresponding relation for $R[X]$-modules, by ([35] Appendix 4 and Proposition 1). □

**Remark 2.9**

In Proposition 2.7 we made crucial use of the fact that

$$F_{\mathbf{Z}_l[G]}(M) = \tau(F_{\mathbf{Z}_l[G]}(Hom(M; \mathbf{Q}_l/\mathbf{Z}_l)))$$

when $M$ is a finite $\mathbf{Z}_l[G]$-module and the abelian group $G$ has cyclic Sylow $l$-subgroup.

Incidentally, Proposition 2.8 seems to contradict ([7] equation (1) p.464)!

**2.10** *Truncation and dualisation*

We shall require one further construction before beginning the proof of Theorem 2.4. Recall from the discussion of §2.1 that both $[\oplus_j F_{2j}, X, \oplus_j F_{2j+1}] \in K_0(\mathbf{Z}_l[G], \mathbf{Q}_l)$ and $det(X) \in \mathbf{Q}_l[G]^*/\mathbf{Z}_l[G]^*$ are independent of the choices of $\mathbf{Q}_l[G]$-module spittings used to construct $X$ and depend only on the quasi-isomorphism class of the complex, $F_*$. Hence we may replace $F_*$ by the quasi-isomorphic truncated complex

$$0 \longrightarrow B_1 \longrightarrow F_1 \longrightarrow F_0 \longrightarrow 0$$

which is perfect. In fact each $\mathbf{Z}_l$-module may be assumed free, by the discussion of §2.5. This may also be seen directly. If we choose the same splittings as in Example 2.2 and §2.5 we obtain two isomorphisms of the form

$$Y : F_0 \otimes \mathbf{Q}_l \oplus B_1 \otimes \mathbf{Q}_l \xrightarrow{\cong} F_1 \otimes \mathbf{Q}_l.$$

$$Z : B_1 \otimes \mathbf{Q}_l \oplus F_3 \otimes \mathbf{Q}_l \oplus F_5 \otimes \mathbf{Q}_l \oplus \ldots \xrightarrow{\cong} F_2 \otimes \mathbf{Q}_l \oplus F_4 \otimes \mathbf{Q}_l \oplus \ldots$$

and $X = (Y \oplus 1 \oplus 1 \oplus \ldots)(1 \oplus Z^{-1})$. Hence $det(X) = det(Y)det(Z)^{-1} \in \mathbf{Q}_l[G]^*$. However, the splittings used to define $Z$ are all $l$-adically integral and so $det(Z) \in \mathbf{Z}_l[G]^*$. Therefore

$$det(X) = det(Y) \in \frac{\mathbf{Q}_l[G]^*}{\mathbf{Z}_l[G]^*},$$

as claimed.

Now set $F^0 = Hom(F_0, \mathbf{Z}_l)$, $F^1 = Hom(F_1, \mathbf{Z}_l)$ and $F^2 = Hom(B_1, \mathbf{Z}_l)$ and consider the cochain complex

$$0 \longrightarrow F^0 \xrightarrow{d_1^*} F^1 \longrightarrow F^2 \longrightarrow 0.$$

From the universal coefficient formula the cohomology groups, $H^i(F^*)$, of this cochain complex vanish except for $i = 1, 2$ in which case

$$H^i(F^*) \cong Ext^1(H_{i-1}(F_*), \mathbf{Z}_l) \cong Hom(H_{i-1}(F_*), \mathbf{Q}_l/\mathbf{Z}_l).$$

The class in $K_0(\mathbf{Z}_l[G], \mathbf{Q}_l)$ defined by $F^*$ depends only on the $[\oplus_j F_{2j}, X, \oplus_j F_{2j+1}]$.

**2.11** *Proof of Theorem 2.4*

Firstly, we observe that the statement depends only on the coset

$$det(X) \in \frac{\mathbf{Q}_l[G]^*}{\mathbf{Z}_l[G]^*}$$

or, equivalently, only on the element $[\oplus_j F_{2j}, X, \oplus_j F_{2j+1}] \in K_0(\mathbf{Z}_l[G], \mathbf{Q}_l)$. Hence, by the discussion of §2.5 and §2.10 we may assume that the $\mathbf{Z}_l[G]$-modules, $F_0, F_1, \ldots, F_k$ and $B_1$ are all free.

We may assume that $t_i \in \mathbf{Q}_p[G]^* \cap \text{ann}_{\mathbf{Z}_i[G]}(H_i(F_*))$. For if $t_i \in \text{ann}_{\mathbf{Z}_i[G]}(H_i(F_*))$ but is not a unit in $\mathbf{Q}_l[G]^*$ we may find a large multiple, $n$, of $|H_0(F_*)| \cdot |H_1(F_*)|$ so that $n + t_i \in \mathbf{Q}_l[G]^*$. Then the result for $t_i^{m_i}$ follows from the that for $(n + t_i)^{m_i}$, by the binomial theorem.

Now consider the result for $t_0 \in \mathbf{Q}_l[G]^* \cap \text{ann}_{\mathbf{Z}_l[G]}(H_0(F_*))$. By Proposition 2.7

$$det(X)t_0^{rank_{\mathbf{Z}_l[G]}(F_0)} \in \text{ann}_{\mathbf{Z}_l[G]}(H_1(F_*)) \lhd \mathbf{Z}_l[G].$$

in part (i) and in part (ii)

$$det(X)t_0^{rank_{\mathbf{Z}_l[G]}(F_0)} \in F_{\mathbf{Z}_l[G]}(H_1(F_*)) \lhd \mathbf{Z}_l[G].$$

Therefore, to complete the proof for $t_0$, it suffices to show that it is possible to choose $F_*$ such that $rank_{\mathbf{Z}_l[G]}(F_0) = m_0$ without changing the class in $K_0(\mathbf{Z}_l[G], \mathbf{Q}_l)$. However, it is easy to construct a perfect chain complex $P_*$, by induction on degree, which is homotopy equivalent to $F_*$ and has $P_0 = \oplus_{i=1}^{m_0} \mathbf{Z}_l[G]$ (for example, see [50] Proposition 2.2.4).

We shall obtain the statement concerning $t_1$ by the same argument applied to the complex
$$0 \longrightarrow F^0 \xrightarrow{d_1^*} F^1 \longrightarrow F^2 \longrightarrow 0$$
which was constructed the truncation and dualisation procedure of §2.10. From the universal coefficient formula the cohomology groups of this cochain complex vanish except for $H^1(F^*) \cong Hom(H_0(F_*), \mathbf{Q}_l/\mathbf{Z}_l)$ and $H^2(F^*) \cong Hom(H_1(F_*), \mathbf{Q}_l/\mathbf{Z}_l)$.

In order to apply the previous argument to this free complex we first observe that, after tensoring with $\mathbf{Q}_l$, the matrix of the dual of
$$Y : F_0 \otimes \mathbf{Q}_l \oplus B_1 \otimes \mathbf{Q}_l \xrightarrow{\cong} F_1 \otimes \mathbf{Q}_l$$
of §2.10, with respect to the dual $\mathbf{Z}_l[G]$-basis, is equal to $\tau(Y^{tr})$ where $\tau$ is the automorphism of the group-ring induced by sending each element of $G$ to its inverse, $\tau(\sum_{g \in G} n_g \cdot g) = \sum_{g \in G} n_g \cdot g^{-1}$, and $Y^{tr}$ is the transpose the matrix of $Y$. Therefore ,the class of $F^*$ in $K_0(\mathbf{Z}_l[G], \mathbf{Q}_l)$ corresponds to $det(\tau(Y^{tr}))^{-1}$, since the dual of $Y$ is an isomorphism of the form
$$F^1 \otimes \mathbf{Q}_l \xrightarrow{\cong} F^0 \otimes \mathbf{Q}_l \oplus F^2 \otimes \mathbf{Q}_l.$$

However, the action of $G$ on $Hom(H_i(F_*), \mathbf{Q}_l/\mathbf{Z}_l)$ is given by $g(f) = f(g^{-1} \cdot -)$ which may be identified with $G$ acting on $H_i(F_*)$ with the new action, $g(z) = \tau(g) \cdot z$. Hence, $t \in \mathbf{Q}_l[G]^* \cap \mathrm{ann}_{\mathbf{Z}_l[G]}(H^2(F^*))$ if and only if $\tau(t) \in \mathbf{Q}_l[G]^* \cap \mathrm{ann}_{\mathbf{Z}_l[G]}(H_1(F_*))$. As before we may assume that $t_1$ lies in $\mathbf{Q}_l[G]^*$ and therefore we may apply Proposition 2.7 to $F^*$ with $t = \tau(t_1)$ to obtain an inclusion of the form
$$det(\tau(Y^{tr}))^{-1}\tau(t_1)^{rank_{\mathbf{Z}_l[G]}(F^2)} \in \mathrm{ann}_{\mathbf{Z}_l[G]}(H^1(F^*)) \triangleleft \mathbf{Z}_l[G]$$
in part (i) and
$$det(\tau(Y^{tr}))^{-1}\tau(t_1)^{rank_{\mathbf{Z}_l[G]}(F^2)} \in F_{\mathbf{Z}_l[G]}(H^1(F^*)) \triangleleft \mathbf{Z}_l[G]$$
in part (ii). Furthermore, by changing $F^*$ as before, we assume that $rank_{\mathbf{Z}_l[G]}(F^2) = m_1$. Finally, since in part (i)
$$\tau(\mathrm{ann}_{\mathbf{Z}_l[G]}(H^1(F^*))) = \mathrm{ann}_{\mathbf{Z}_l[G]}(H_0(F_*))$$
and in part (ii)
$$\tau(F_{\mathbf{Z}_l[G]}(H^1(F^*))) = F_{\mathbf{Z}_l[G]}(H_0(F_*))$$
we may apply $\tau$ to this inclusion to obtain, as required,
$$det(X)^{-1} t_1^{m_1} \in \mathrm{ann}_{\mathbf{Z}_l[G]}(H_0(F_*)) \triangleleft \mathbf{Z}_l[G]$$

in part (i) and in part (ii)

$$det(X)^{-1}t_1^{m_1} \in F_{\mathbf{Z}_l[G]}(H_0(F_*)) \triangleleft \mathbf{Z}_l[G],$$

since $det(\tau(Y^{tr})) = det(\tau(Y))$ differs from $\tau(det(X))$ by a unit in the $l$-adic group-ring. □.

**Corollary 2.12**

Suppose, in the situation of Theorem 2.4, that $m_1 = 1$ (that is, the $\mathbf{Z}_l[G]$-module $Hom(H_1(F_*), \mathbf{Q}_l/\mathbf{Z}_l)$ is generated by one element) then:

(i)

$$\{t^{m_0} \mid t \in \operatorname{ann}_{\mathbf{Z}_l[G]}(H_0(F_*))\} \subseteq det(X)^{-1} \cdot \operatorname{ann}_{\mathbf{Z}_l[G]}(H_1(F_*)) \subseteq \operatorname{ann}_{\mathbf{Z}_l[G]}(H_0(F_*))$$

(ii) if the Sylow $l$-subgroup of $G$ is cyclic then in (i) the right-hand $\operatorname{ann}_{\mathbf{Z}_l[G]}(H_0(F_*))$ may be replaced by $F_{\mathbf{Z}_l[G]}(H_0(F_*))$.

# 3 Function fields

**3.1** Suppose that $N/K$ is a finite Galois extension of global function fields with abelian Galois group $G(N/K)$ and let $p = char(N)$. Let $S$ be a finite, non-empty set of primes of $K$, including all those primes which ramify in $N/K$. Let $S'$ denote the set of primes of $N$ lying over those of $S$. Let $O_{N,S'}$ be the ring of $S'$-integers of $N$; that is, the elements of $N$ which are regular away from $S'$. Then $W = \operatorname{Spec}(O_{N,S'})$ is a smooth affine curve over $\mathbf{F}_p$.

Let $l$ be a prime different from $p$. If $\overline{\mathbf{Q}_l}$ denotes an algebraic closure of $\mathbf{Q}_l$ and $\chi : G(L/K) \longrightarrow \overline{\mathbf{Q}_l}^*$ is a one-dimensional $l$-adic representation we denote by $L_{K,S}(t,\chi)$ denote the Artin L-function ([36] p.293; see also [34]) with the Euler factors associated to primes of $S$ removed. In the notation of [36] $L_K(t,\chi)$ would be denoted by $L(W, \chi, t)$.

**Proposition 3.2**

Let $N/K$ be as in §3.1. Then, for each integer $n \geq 2$, there exists a (unique) Stickelberger element

$$\Theta_{N/K,S}(n) \in \mathbf{Q}_l[G(N/K)]^*$$

such that

$$\chi(\Theta_{N/K,S}(n)) = L_{K,S}(1-n, \chi^{-1})$$

for all one-dimensional $l$-adic representations, $\chi$.

**Proof**

In ([8] Proposition 7.13) it is shown that there exists an element

$$\Omega_{n-1,l}(N/K) \in K_0(\mathbf{Z}_l[G], \mathbf{Q}_l) \cong \frac{\mathbf{Q}_l[G]^*}{\mathbf{Z}_l[G]^*}$$

which corresponds to $\chi \mapsto L_{K,S}(1-n, \chi^{-1})^{-1}$. This means that there exists a unit, $\Omega \in \mathbf{Q}_l[G]^*$, such that $\chi(\Omega) = L_{K,S}(1-n, \chi^{-1})^{-1}$ for all $\chi$. Therefore we set $\Theta_{N/K,S}(n) = \Omega^{-1}$. □

**3.3** For $i = 0, 1$ and $n \geq 1$ we have étale cohomology groups ([8] §7)

$$H^i_{\text{ét}}(W; \mathbf{Q}_l/\mathbf{Z}_l(n+1)) \cong H^{i+1}_{\text{ét}}(W; \mathbf{Z}_l(n+1))$$

where $W = \text{Spec}(O_{N,S'})$, as in §3.1. Each of these finite groups is a $\mathbf{Z}_l[G(N/K)]$-module.

Suppose that $E/K$ is a Galois extension containing $N$ and having $\mu_{l^\infty}(E)$, the group of $l$-primary roots of unity in $E$, very large. Then

$$H^0_{\text{ét}}(W; \mathbf{Q}_l/\mathbf{Z}_l(n+1)) \cong (\mu_{l^\infty}(E)^{\otimes n+1})^{G(E/N)},$$

which is a cyclic group.

The annihilator, $\text{ann}_{\mathbf{Z}_l[G(N/K)]}(H^0_{\text{ét}}(W; \mathbf{Q}_l/\mathbf{Z}_l(n+1)))$ is well-known, being computed as in ([10] Lemma 2.3; [54] p.82). This annihilator is related to $\text{ann}_{\mathbf{Z}_l[G(N/K)]}(H^1_{\text{ét}}(W; \mathbf{Q}_l/\mathbf{Z}_l(n+1)))$ and the Fitting ideal $F_{\mathbf{Z}_l[G(N/K)]}(H^1_{\text{ét}}(W; \mathbf{Q}_l/\mathbf{Z}_l(n+1)))$ by the following result.

**Theorem 3.4**

*Suppose that $N/K$ is a finite Galois extension of global function fields with abelian Galois group, $G(N/K)$. If $n \geq 1$ with $l$ and $W = Spec(O_{N,S'})$ as in §3.1 then*

(i)

$$\{t^m \mid t \in \text{ann}_{\mathbf{Z}_l[G(N/K)]}(H^1_{\text{ét}}(W; \mathbf{Q}_l/\mathbf{Z}_l(n+1)))\}$$

$$\subseteq \Theta_{N/K,S}(n+1) \cdot \text{ann}_{\mathbf{Z}_l[G(N/K)]}(H^0_{\text{ét}}(W; \mathbf{Q}_l/\mathbf{Z}_l(n+1)))$$

$$\subseteq \text{ann}_{\mathbf{Z}_l[G(N/K)]}(H^1_{\text{ét}}(W; \mathbf{Q}_l/\mathbf{Z}_l(n+1)))$$

*where $m$ is the minimum number of generators for the $\mathbf{Z}_l[G(N/K)]$-module, $H^1_{\text{ét}}(W; \mathbf{Q}_l/\mathbf{Z}_l(n+1))$.*

*(ii) If the Sylow $l$-subgroup of $G(N/K)$ is cyclic then the final $\text{ann}_{\mathbf{Z}_l[G(N/K)]}(H^1_{\text{ét}}(W; \mathbf{Q}_l/\mathbf{Z}_l(n+1)))$ in (i) may be replaced by the Fitting ideal $F_{\mathbf{Z}_l[G(N/K)]}(H^1_{\text{ét}}(W; \mathbf{Q}_l/\mathbf{Z}_l(n+1)))$.*

**Proof**

In ([8] Proposition 7.7) it is shown that there exists a bounded, perfect cochain complex of $\mathbf{Z}_l[G(N/K)]$-modules of the form

$$M^* : 0 \to M^s \to M^{s+1} \to \ldots \to M^t \to 0$$

together with a $\mathbf{Z}_l[G(N/K)]$-module automorphism, $\alpha : M^* \to M^*$, such that for each $j$ the homomorphism $1 - \alpha : M^j \to M^j$ is injective with finite cokernel. In addition, the mapping cone, $C^* = Cone(1 - \alpha)$ has trivial cohomology with the exception of $H^i(C^*) \cong H^i_{\text{ét}}(X; \mathbf{Q}_l/\mathbf{Z}_l(n+1))$ for $i = 0, 1$. By the construction of Example 2.2 the complex $C^*$ gives rise to an element of

$$K_0(\mathbf{Z}_l[G], \mathbf{Q}_l) \cong \frac{\mathbf{Q}_l[G]^*}{\mathbf{Z}_l[G]^*}.$$

In fact the truncated complex, $D^*$, given by $D^j = C^j$ for $j \leq 0$, $D^1 = Ker(d : C^1 \to C^2)$ and $D^j = 0$ for $j \geq 2$ defines the same element. Now setting $F_{1-j} = D^j$ and applying Theorem 2.4 and Corollary 2.12 to $F_*$ yields, in case (i),

$$\{z^t \mid z \in \text{ann}_{\mathbf{Z}_l[G(N/K)]}(H^1_{\text{ét}}(X; \mathbf{Q}_l/\mathbf{Z}_l(n+1)))\}$$

$$\subseteq det(X)^{-1} \cdot \text{ann}_{\mathbf{Z}_l[G(N/K)]}(H^0_{\text{ét}}(X; \mathbf{Q}_l/\mathbf{Z}_l(n+1)))$$

$$\subseteq \text{ann}_{\mathbf{Z}_l[G(N/K)]}(H^1_{\text{ét}}(X; \mathbf{Q}_l/\mathbf{Z}_l(n+1)))$$

where $[\oplus_j F_{2j}, X, \oplus_j F_{2j+1}]$ is the element of $K_0(\mathbf{Z}_l[G], \mathbf{Q}_l)$ given by the construction of Example 2.2. However, for the cochain complex $C^*$ it is shown in ([8] Proposition 7.17) that this element satisfies $\chi(det(X)) = L_{K,S}(-n, \chi^{-1})^{-1}$ for all one-dimensional $l$-adic representations, $\chi$. Hence

$$det(X) = \Theta_{N/K,S}(n+1)^{-1},$$

as required in (i). In case (ii) we are permitted by Corollary 2.12(ii) to replace the final annihilator by the Fitting ideal. □

**Corollary 3.5**

*Let $\sqrt{I}$ denote the radical of an ideal $I$. Then, in Theorem 3.4,*

$$\sqrt{\Theta_{N/K,S}(n+1) \cdot \text{ann}_{\mathbf{Z}_l[G(N/K)]}(H^0_{\text{ét}}(W; \mathbf{Q}_l/\mathbf{Z}_l(n+1)))}$$

$$= \sqrt{\text{ann}_{\mathbf{Z}_l[G(N/K)]}(H^1_{\text{ét}}(W; \mathbf{Q}_l/\mathbf{Z}_l(n+1)))}.$$

**Remark 3.6** (i) When $l$ is an odd prime, Corollary 3.5 is also true in the case of a Galois extension, $L/K$, of number fields with abelian Galois group but is much harder to prove. I learnt the proof given in ([50] Theorem 7.3.2) from Ted Chinburg.

(ii) Let $S_1$ be any finite, non-empty set of primes of $K$ and let $S_1'$ denote the set of primes of $N$ lying over those of $S_1$. If $S_1 \subseteq S$ with $S$ as in Theorem 3.4 then the localisation sequence in étale cohomology shows that

$$\Theta_{N/K,S}(n+1) \cdot ann_{\mathbf{Z}_l[G(N/K)]}(H^0_{\text{ét}}(W; \mathbf{Q}_l/\mathbf{Z}_l(n+1)))$$

$$\subseteq ann_{\mathbf{Z}_l[G(N/K)]}(H^1_{\text{ét}}(\text{Spec}(\mathcal{O}_{N,S_1'}); \mathbf{Q}_l/\mathbf{Z}_l(n+1))).$$

This raises the more subtle question of whether one can replace $\Theta_{N/K,S}(n+1)$ by $\Theta_{N/K,S_1}(n+1)$ to give an inclusion similar to that of Conjecture 1.2; that is, not involving any ramification condition on the set of primes, $S_1$. Of course, often the answer to this question is affirmative for trivial reasons.

(iii) As in the number field case mentioned in Remark 1.4, there is an equivariant Chern class homomorphism of the form [14]

$$c_{n+1,2}: K_{2n}(W) \otimes \mathbf{Z}_l \longrightarrow H^1_{\text{ét}}(W; \mathbf{Q}_l/\mathbf{Z}_l(n+1)),$$

which the Lichtenbaum-Quillen Conjecture predicts to be an isomorphism. It is known that $c_{1,2}$ is an isomorphism.

# 4 Absolutely abelian number fields

**4.1** *The Determinant Functor*

We quickly review some formalism of the determinant functor as developed in [24].

For any commutative ring $R$, which will be $\mathbf{Z}_l[G]$ and $\mathbf{Q}_l[G]$ in our applications, a graded invertible module is a pair, $(L, \alpha)$, consisting of an invertible (i.e. projective of rank one) $R$-module $L$ and a locally constant function $\alpha: Spec(R) \to \mathbf{Z}$. A morphism $h: (L, \alpha) \to (M, \beta)$ is an $R$-module homomorphism, $h: L \to M$, such that $\alpha(P) \neq \beta(P)$ implies $h(P) = 0$ for all $P \in Spec(R)$. Let $\mathcal{P}_R$ denote the category of graded invertible modules and their isomorphisms. Then $\mathcal{P}_R$ is a symmetric monoidal category in the sense of [33] with tensor product defined by $(L, \alpha) \otimes (M, \beta) = (L \otimes_R M, \alpha + \beta)$, the usual associativity constraint, unit object $(R, 0)$ and commutativity constraint

$$\psi(l \otimes m) := \psi_{(L,\alpha),(M,\beta)}(l \otimes m) = (-1)^{\alpha(P)\beta(P)} m \otimes l$$

for local sections $l \in L_P, m \in M_P$. We define $(L, \alpha)^{-1} := (Hom_R(L, R), -\alpha)$. This definition extends to a covariant tensor functor $\mathcal{P}_R \to \mathcal{P}_R$ by setting $h^{-1} := Hom_R(h, R)^{-1}$.

For a finitely generated projective $R$-module, $M$, we define an object in $\mathcal{P}_R$ by

$$D_R(M) = (det_R M, rank_R(M))$$

where $det_R M = \bigwedge_R^{rank}(M)$ denotes the maximal exterior power. For a bounded perfect complex $M^*$ we set

$$D_R(M^*) = \bigotimes_{i \in \mathbf{Z}} D_R(M^i)^{(-1)^{i+1}},$$

the tensor product of graded invertible $R$-modules. Here we use the parity convention of [5] rather than that of [24].

We write $\mathcal{D}(R)$ for the derived category of the homotopy category of bounded complexes of $R$-modules, and $\mathcal{D}^p(R)$ for the full triangulated subcategory of perfect complexes of $R$-modules. We write $\mathcal{D}^{pis}(R)$ for the subcategory of $\mathcal{D}^p(R)$ with the same objects but the morphisms are restricted to be quasi-isomorphisms. If $R$ is reduced, then $D_R$ extends to a functor from $\mathcal{D}^{pis}(R)$ to $\mathcal{P}_R$ such that for every distinguished triangle $C_1 \to C_2 \to C_3 \to C_1[1]$ in $\mathcal{D}^p(R)$ there is an isomorphism of the form

$$D_R(C_1) \otimes D_R(C_3) \stackrel{\cong}{\to} D_R(C_2)$$

which is functorial in the triangle ([24] Proposition 7).

If $C$ is an acyclic perfect complex then there exists a canonical isomorphism of the form

$$D_R(C) \stackrel{\cong}{\Longrightarrow} R.$$

For example, suppose that $G$ is a finite abelian group and that $F^{-1}, F^0$ are two finitely generated, free $\mathbf{Z}_l[G]$-modules and that

$$0 \to F^{-1} \stackrel{Z}{\to} F^0 \to 0$$

is a complex of $\mathbf{Z}_l[G]$-modules which becomes exact when tensored with $\mathbf{Q}_l$. Then, by choosing $\mathbf{Z}_l[G]$-bases we obtain a well-defined element

$$det(Z) \in \frac{\mathbf{Q}_l[G]^*}{\mathbf{Z}_l[G]^*}.$$

On the other hand, if $rank_{\mathbf{Z}_l[G]} F^{-1} = m$, we have an exact sequence of rank one modules of the form

$$0 \to \bigwedge^m F^{-1} \stackrel{det(Z)}{\to} \bigwedge^m F^0$$

or, equivalently,

$$0 \to D_{\mathbf{Z}_l[G]}(F^{-1}) \otimes D_{\mathbf{Z}_l[G]}(F^0)^{-1} \stackrel{det(Z)}{\to} \mathbf{Z}_l[G].$$

The image of this homomorphism under the canonical isomorphism

$$D_{\mathbf{Q}_l[G]}(F^{-1}) \otimes D_{\mathbf{Q}_l[G]}(F^0)^{-1} \cong \mathbf{Q}_l[G]$$

is, by definition, the invertible module $D_{\mathbf{Z}_l[G]}(F^*)$ which is therefore the free $\mathbf{Z}_l[G]$-submodule of $\mathbf{Q}_l[G]$ generated by $det(Z)$. In the notation of Example 2.2, $det(Z) = det(X)^{-1}$. In general for any bounded complex of finitely generated, free $\mathbf{Z}_l[G]$-modules with finite cohomology groups, as in Example 2.2, we have
$$D_{\mathbf{Z}_l[G]}(F^*) = \mathbf{Z}_l[G] < det(X)^{-1} > \subset \mathbf{Q}_l[G].$$

**4.2** We recall now some crucial results from [5]. Let $l$ be an odd prime, $m$ a positive integer not divisible by $l$ and $r$ a strictly negative integer. If $\xi_t = e^{2\pi i/t}$ we have a canonical projection of the form $G(\mathbf{Q}(\xi_{ml^{s+1}})/\mathbf{Q}) \longrightarrow G(\mathbf{Q}(\xi_{ml^s})/\mathbf{Q})$ and, taking the inverse limit over the induced homomorphisms of $l$-adic group-rings, we define
$$\Lambda_m^\infty = \lim_{\substack{\longleftarrow \\ n}} \mathbf{Z}_l[G(\mathbf{Q}(\xi_{ml^{n+1}})/\mathbf{Q})].$$

Write $Q(\Lambda_m^\infty)$ for the total quotient ring of $\Lambda_m^\infty$ ([16] p.60). In ([5] §7 (48)) a bounded perfect complex of $\Lambda_m^\infty$-modules, denoted by $C_m(r)$, is constructed for which ([5] Lemma 7.2) $C_m(r) \otimes_{\Lambda_m^\infty} Q(\Lambda_m^\infty)$ is acyclic. Therefore we have a canonical isomorphism of the form
$$D_{\Lambda_m^\infty}(C_m(r)) \otimes_{\Lambda_m^\infty} Q(\Lambda_m^\infty) \xrightarrow{\cong} Q(\Lambda_m^\infty).$$

We shall recall the definition of $C_m(r)$ presently. Furthermore ([5] Theorem 7.1) calculates a $\Lambda_m^\infty$-basis for the free, rank one $\Lambda_m^\infty$-module $D_{\Lambda_m^\infty}(C_m(r))$ by showing that the imagine of $D_{\Lambda_m^\infty}(C_m(r))$ under the canonical homomorphism into $Q(\Lambda_m^\infty)$ is $\Lambda_m\langle(e_+ + e_- \cdot g_m)\rangle$ in the notation of [5]. Fix an integer $s \geq 0$ and a topological generator $\gamma_s \in G(\mathbf{Q}(\xi_{ml^\infty}/\mathbf{Q}(\xi_{ml^{s+1}}))$ and form the bounded perfect complex of $\mathbf{Z}_l[G(\mathbf{Q}(\xi_{ml^{s+1}})/\mathbf{Q})]$-modules
$$C_{0,m}(r) = C(r) \otimes_{\Lambda_m^\infty}^{\mathbf{L}} \mathbf{Z}_l[G(\mathbf{Q}(\xi_{ml^{s+1}})/\mathbf{Q})].$$

so that
$$C_m(r) \xrightarrow{1-\gamma_s} C_m(r) \longrightarrow C_{0,m}(r) \longrightarrow C_m(r)[1]$$

is a distinguished triangle in the derived category of the homotopy category of $\Lambda_m^\infty$-modules. Since $l$ is odd, taking coinvariants (or equivalently invariants) under complex conjugation $C_{0,m}(r) \mapsto C_{0,m}(r)^+$ is exact and results in a bounded perfect complex of $\mathbf{Z}_l[G(\mathbf{Q}(\xi_{ml^{s+1}})^+/\mathbf{Q})]$-modules where $\mathbf{Q}(\xi_{ml^{s+1}})^+$ is the totally real subfield of $\mathbf{Q}(\xi_{ml^{s+1}})$. In addition, the cohomology groups of $C_{0,m}(r)^+$ are finite so that we have a canonical isomorphism

$$D_{\mathbf{Z}_l[G(\mathbf{Q}(\xi_{ml^{s+1}})^+/\mathbf{Q})]}(C_{0,m}(r)^+) \otimes_{\mathbf{Z}_l[G(\mathbf{Q}(\xi_{ml^{s+1}})^+/\mathbf{Q})]} \mathbf{Q}_l[G(\mathbf{Q}(\xi_{ml^{s+1}})^+/\mathbf{Q})]$$

$$\xrightarrow{\cong} \mathbf{Q}_l[G(\mathbf{Q}(\xi_{ml^{s+1}})^+/\mathbf{Q})].$$

Let $\tau$ denote complex conjugation and let
$$\chi : G(\mathbf{Q}(\xi_{ml^\infty})/\mathbf{Q}) \longrightarrow G(\mathbf{Q}(\xi_{ml^{s+1}})^+/\mathbf{Q}) \longrightarrow \overline{\mathbf{Q}}_l^*$$
denote an $l$-adic character. Suppose that $z_{m,r}$ is a basis element for the free rank one $\mathbf{Z}_l[G(\mathbf{Q}(\xi_{ml^{s+1}})^+/\mathbf{Q})]$-module
$$D_{\mathbf{Z}_l[G(\mathbf{Q}(\xi_{ml^{s+1}})^+/\mathbf{Q})]}(C_{0,m}(r)^+) \subset \mathbf{Q}_l[G(\mathbf{Q}(\xi_{ml^{s+1}})^+/\mathbf{Q})].$$
By the descent lemma ([5] Lemma 8.1) with $C$ replaced by $C_m(r)$) $z_{m,r}$ may be chosen to satisfy
$$\chi(z_{m,r}) = \chi((1+(-1)^r\tau)/2) + \chi((1+(-1)^{r-1}\tau)/2)\chi(g_m(r))$$
$$= \begin{cases} 1 & \text{if } r \text{ is even,} \\ -L(r,\tilde{\chi}) & \text{if } r \text{ is odd} \end{cases}$$
since $-\chi(g_m(r)) = L(r,\tilde{\chi})$, by ([5] Proposition 5.5), the Dirichlet L-function at $s = r$ of the primitive Dirichlet character $\tilde{\chi}$ associated to $\chi$.

Let $X_l = \mathrm{Spec}(\mathbf{Z}[\xi_{ml^{s+1}}][1/ml])$, $X_l^+ = \mathrm{Spec}(\mathbf{Z}[\xi_{ml^{s+1}}]^+[1/ml])$ denote the spectrum of the $S_{ml}$-integers of $\mathbf{Q}(\xi_{ml^{s+1}})$, $\mathbf{Q}(\xi_{ml^{s+1}})^+$, respectively. Here, as usual, $S_{ml}$ equals the set of primes dividing $ml$.

Now we must examine the construction of $C_m(r)$ prior to applying Theorem 2.4 to $C_{0,m}(r)^+$. Following [5] let $\Sigma_{L_m^\infty}$ denote the set of embeddings of $\mathbf{Q}(\xi_{ml^\infty})$ into the complex numbers. Then complex conjugation acts diagonally on $\prod_{\Sigma_{L_m^\infty}} \mathbf{Z}_l(-r)$ so that the fixed points $(\prod_{\Sigma_{L_m^\infty}} \mathbf{Z}_l(-r))^+$ naturally form a projective $\Lambda_m^\infty$-module generated by $e_+(\tilde{\sigma}_m^\infty \nu_r)$ in the notation of ([5] §7). In ([5] §7 (48)) a $\Lambda_m^\infty$-module homomorphism $c_m(r)$ from $(\prod_{\Sigma_{L_m^\infty}} \mathbf{Z}_l(-r))^+$ to $H^1_{\text{ét}}(X_l; \mathbf{Z}_l(1-r))$ sending $e_+(\tilde{\sigma}_m^\infty \nu_r)$ to a "cyclotomic element" denoted by $\eta_m(r)$. This uniquely determines a homomorphism in the derived category
$$c_m(r) \in Hom_{\mathcal{D}(\Lambda_m^\infty)}((\prod_{\Sigma_{L_m^\infty}} \mathbf{Z}_l(-r))^+[-1], R\Gamma(X_l; \mathbf{Z}_l(1-r)))$$
and $C_m(r)$ is the bounded perfect complex of $\Lambda_m^\infty$-modules given by the mapping cone of $c_m(r)$. Hence there is a distinguished triangle of perfect complexes in $\mathcal{D}(\Lambda_m^\infty)$ of the form
$$(\prod_{\Sigma_{L_m^\infty}} \mathbf{Z}_l(-r))^+[-1] \xrightarrow{c_m(r)} R\Gamma(X_l; \mathbf{Z}_l(1-r)) \longrightarrow C_m(r).$$
The complex $C_{0,m}(r)$ is obtained by applying $(- \otimes_\Lambda^{\mathrm{L}} \mathbf{Z}_l[G(\mathbf{Q}(\xi_{ml^{s+1}})/\mathbf{Q})])$. Complex conjugation acts on $\mathbf{Z}_l(-r)$ as multiplication by $(-1)^r$ so that, when $r$ is odd, there are $\mathbf{Z}_l[G(\mathbf{Q}(\xi_{ml^{s+1}})^+/\mathbf{Q})]$-module isomorphisms
$$H^i(C_{0,m}(r)^+) \cong H^i_{\text{ét}}(X_l^+; \mathbf{Z}_l(1-r)) \cong H^{i-1}_{\text{ét}}(X_l^+; (\mathbf{Q}_l/\mathbf{Z}_l)(1-r))$$

for $i = 1, 2$ and $H^i(C_{0,m}(r)^+) = 0$ otherwise. However, when $r$ is even, the action by complex conjugation is trivial and the cohomology long exact sequence yields a short exact sequence

$$0 \longrightarrow \mathbf{Z}_l[G(\mathbf{Q}(\xi_{ml^{s+1}})^+/\mathbf{Q})]\langle e_+(\tilde{\sigma}_m^\infty \nu_r)\rangle \overset{c_m(r)}{\longrightarrow} H^1_{\text{ét}}(X_l^+; \mathbf{Z}_l(1-r))$$

$$\longrightarrow H^1(C_{0,m}(r)^+) \longrightarrow 0$$

and an isomorphism

$$H^2_{\text{ét}}(X_l^+; \mathbf{Z}_l(1-r)) \overset{\cong}{\longrightarrow} H^2(C_{0,m}(r)^+)$$

while $H^i(C_{0,m}(r)^+) = 0$ otherwise.

**4.3** *Cyclotomic fields – proof of Theorem 1.6*

The discussion of §4.2 shows that $C_{0,m}(r)^+$ is a perfect complex of $\mathbf{Z}_l[G(\mathbf{Q}(\xi_{ml^{s+1}})^+/\mathbf{Q})]$ for which, when $r = -1, -3, -5, \ldots$,

$$det(X)^{-1} = -\Theta_{\mathbf{Q}(\xi_{ml^{s+1}})^+/\mathbf{Q}}(1-r) \in \frac{\mathbf{Q}_l[G(\mathbf{Q}(\xi_{ml^{s+1}})^+/\mathbf{Q})]^*}{\mathbf{Z}_l[G(\mathbf{Q}(\xi_{ml^{s+1}})^+/\mathbf{Q})]^*}.$$

Since the non-zero (finite) homology groups of $C_{0,m}(r)^+$ are given by $H_i = H^{1-i}_{\text{ét}}(X_l^+; (\mathbf{Q}_l/\mathbf{Z}_l)(1-r))$ for $i = 0, 1$ and $H_1$ is a cyclic group the result follows from Theorem 2.4 and Corollary 2.12. □

**Remark 4.4** *Base change for higher Stickelberger ideals*

The main theorem of [45] together with Corollary 1.7 shows that Corollary 1.7 remains true if $\mathbf{Q}(\xi_{ml^{s+1}})^+/\mathbf{Q}$ is replaced by any intermediate extension $\mathbf{Q}(\xi_{ml^{s+1}})^+/K$. For another explanation of base change phenomena see ([50] Remark 7.2.11(iii)).

**4.5** $H^i_{\text{ét}}(X_l^+; \mathbf{Z}_l(1-r))$, $0 > r$ even

Let $l$ be an odd prime. The discussion of §4.2 serves equally well when applying Theorem 2.4 and Corollary 2.12 to $C_{0,m}(r)^+$ with $r$ even but the outcome it slightly different. Since complex conjugation $\tau$ acts as multiplication by $-1$ on the finite cyclic $l$-group $H^0_{\text{ét}}(X_l; (\mathbf{Q}_l/\mathbf{Z}_l)(1-r))$ we have $H^0_{\text{ét}}(X_l^+; (\mathbf{Q}_l/\mathbf{Z}_l)(1-r)) = H^0_{\text{ét}}(X_l; (\mathbf{Q}_l/\mathbf{Z}_l)(1-r))^{\langle\tau\rangle} = 0$, as $l$ is odd. Hence $H^1_{\text{ét}}(X_l^+; \mathbf{Z}_l(1-r))$ is $\mathbf{Z}_l$-torsion free because embeds into $H^1_{\text{ét}}(X_l^+; \mathbf{Q}_l(1-r))$. The image of the homomorphism $c_m(r)$ of §4.2 is a free rank one $\mathbf{Z}_l[G(\mathbf{Q}(\xi_{ml^{s+1}})^+/\mathbf{Q})]$-module generated by the "cyclotomic element" $c_m(r)(e_+(\tilde{\sigma}_m^\infty \nu_r)) = \eta_m(r)$ ([5] §7 and §8). Following ([25] §5; see also [37] Part II §1) denote the image of $c_m(r)$ by

$$im(c_m(r)) = H^1_{\text{cyclo}}(X_l^+; \mathbf{Z}_l(1-r)) \subseteq H^1_{\text{ét}}(X_l^+; \mathbf{Z}_l(1-r)),$$

the submodule of cyclotomic elements, then

$$\mathcal{H}^1_{\text{cyclo}} = \frac{H^1_{\text{ét}}(X_l^+; \mathbf{Z}_l(1-r))}{H^1_{\text{cyclo}}(X_l^+; \mathbf{Z}_l(1-r))}$$

is a finite $l$-group ([5] §8). The short exact sequence of $\mathbf{Z}_l[G(\mathbf{Q}(\xi_{ml^{s+1}})^+/\mathbf{Q})]$-modules

$$0 \longrightarrow Hom_{\mathbf{Z}_l}(H^1_{\text{ét}}(X_l^+; \mathbf{Z}_l(1-r)), \mathbf{Z}_l) \longrightarrow$$

$$Hom_{\mathbf{Z}_l}(\mathbf{Z}_l[G(\mathbf{Q}(\xi_{ml^{s+1}})^+/\mathbf{Q})]\langle e_+(\tilde{\sigma}_m^\infty \nu_r)\rangle, \mathbf{Z}_l) \longrightarrow Ext^1_{\mathbf{Z}_l}(\mathcal{H}^1_{\text{cyclo}}, \mathbf{Z}_l) \longrightarrow 0$$

shows that

$$Ext^1_{\mathbf{Z}_l}(\mathcal{H}^1_{\text{cyclo}}, \mathbf{Z}_l) \cong Hom(\mathcal{H}^1_{\text{cyclo}}, \mathbf{Q}_l/\mathbf{Z}_l)$$

is generated by one element, because $\mathbf{Z}_l[G(\mathbf{Q}(\xi_{ml^{s+1}})^+/\mathbf{Q})]$ is a self-dual $\mathbf{Z}_l[G(\mathbf{Q}(\xi_{ml^{s+1}})^+/\mathbf{Q})]$-module.

Therefore we may apply Corollary 2.12 to $C_{0,m}(r)^+$ when $r < 0$ is even to obtain the following result.

**Theorem 4.6**

Let $l$ be an odd prime. Let $l, m$ and $X_l^+$ be as in the notation of §1.5. Let $H^1_{\text{cyclo}}(X_l^+; \mathbf{Z}_l(1-r))$ denote the cyclotomic elements of §4.5. Then, for $r = -2, -4, -6, \ldots$ and any positive integer $s$:

(i) there exists a chain of annihilator ideal relations for étale cohomology of the form

$$\{t^{m_0} \mid t \in \text{ann}_{\mathbf{Z}_l[G(\mathbf{Q}(\xi_{ml^{s+1}})^+/\mathbf{Q})]}(H^2_{\text{ét}}(X_l^+; \mathbf{Z}_l(1-r)))\}$$

$$\subseteq \text{ann}_{\mathbf{Z}_l[G(\mathbf{Q}(\xi_{ml^{s+1}})^+/\mathbf{Q})]}\left(\frac{H^1_{\text{ét}}(X_l^+; \mathbf{Z}_l(1-r))}{H^1_{\text{cyclo}}(X_l^+; \mathbf{Z}_l(1-r))}\right)$$

$$\subseteq \text{ann}_{\mathbf{Z}_l[G(\mathbf{Q}(\xi_{ml^{s+1}})^+/\mathbf{Q})]}(H^2_{\text{ét}}(X_l^+; \mathbf{Z}_l(1-r))).$$

(ii) If $l$ does not divide $m - 1$ then in (i) the final

$$\text{ann}_{\mathbf{Z}_l[G(\mathbf{Q}(\xi_{ml^{s+1}})^+/\mathbf{Q})]}(H^2_{\text{ét}}(X_l^+; \mathbf{Z}_l(1-r)))$$

may be replaced by the Fitting ideal

$$F_{\mathbf{Z}_l[G(\mathbf{Q}(\xi_{ml^{s+1}})^+/\mathbf{Q})]}(H^2_{\text{ét}}(X_l^+; \mathbf{Z}_l(1-r))).$$

Here $m_0$ is the minimal number of generators of the $\mathbf{Z}_l[G(\mathbf{Q}(\xi_{ml^{s+1}})^+/\mathbf{Q})]$-module $H^2_{\text{ét}}(X_l^+; \mathbf{Z}_l(1-r))$.

**Remark 4.7** (i) There is no Stickelberger element in the statement of Theorem 4.6 because the discussion of §4.2 shows that the element of

$$K_0(\mathbf{Z}_l[G(\mathbf{Q}(\xi_{ml^{s+1}})^+/\mathbf{Q})], \mathbf{Q}_l) \cong \frac{\mathbf{Q}_l[G(\mathbf{Q}(\xi_{ml^{s+1}})^+/\mathbf{Q})]^*}{\mathbf{Z}_l[G(\mathbf{Q}(\xi_{ml^{s+1}})^+/\mathbf{Q})]^*}$$

defined by $C_{0,m}(r)^+$ is trivial when $r < 0$ is even (i.e. $det(X) = 1$ in the notation of Theorem 2.4).

Observe that this triviality immediately implies the central result of ([25] Theorem 5.4 and Corollary 5.5). Namely, if $\chi$ is a character of $G(\mathbf{Q}(\xi_{ml^{s+1}})^+/\mathbf{Q})$ of order prime to $l$ then the $\chi$-eigenspaces satisfy

$$[H^1_{\text{ét}}(X_l^+; \mathbf{Z}_l(1-r))^\chi : H^1_{\text{cyclo}}(X_l^+; \mathbf{Z}_l(1-r))^\chi] = |H^2_{\text{ét}}(X_l^+; \mathbf{Z}_l(1-r))^\chi|. \quad (4.7.1)$$

To see this, firstly map

$$0 = [C_{0,m}(r)^+_{\text{even}}, X, C_{0,m}(r)^+_{\text{odd}}] \in K_0(\mathbf{Z}_l[G(\mathbf{Q}(\xi_{ml^{s+1}})^+/\mathbf{Q})], \mathbf{Q}_l)$$

to $G_0(\mathbf{Z}_l[G(\mathbf{Q}(\xi_{ml^{s+1}})^+/\mathbf{Q})], \mathbf{Q}_l)$ by the forgetful map. However, there is an isomorphism

$$G_0(\mathbf{Z}_l[G(\mathbf{Q}(\xi_{ml^{s+1}})^+/\mathbf{Q})], \mathbf{Q}_l) \cong G_0T(\mathbf{Z}_l[G(\mathbf{Q}(\xi_{ml^{s+1}})^+/\mathbf{Q})]),$$

the Grothendieck group of finite $\mathbf{Z}_l[G(\mathbf{Q}(\xi_{ml^{s+1}})^+/\mathbf{Q})]$-modules. Under the decomposition $\mathbf{Z}_l[G] \cong \prod_\chi \mathbf{Z}_l[\chi][G_1]$ of Proposition 2.8(proof) the $\chi$-component yields

$$[\frac{H^1_{\text{ét}}(X_l^+; \mathbf{Z}_l(1-r))^\chi}{H^1_{\text{cyclo}}(X_l^+; \mathbf{Z}_l(1-r))^\chi}] = |H^2_{\text{ét}}(X_l^+; \mathbf{Z}_l(1-r))^\chi|] \in G_0(\mathbf{Z}_l[\chi][G_1])$$

and mapping to $G_0(\mathbf{Z}_l[\chi])$ by the augmentation yields (4.7.1).

(ii) Recall that $\Theta_{\mathbf{Q}(\xi_{ml^{s+1}})^+/\mathbf{Q}}(1-r)$ is zero when $0 > r$ is even [58] so that a result like Theorem 1.6 would not make sense. Next we examine how to manufacture elements of

$$\text{ann}_{\mathbf{Z}_l[G(\mathbf{Q}(\xi_{ml^{s+1}})^+/\mathbf{Q})]}(\frac{H^1_{\text{ét}}(X_l^+; \mathbf{Z}_l(1-r))}{H^1_{\text{cyclo}}(X_l^+; \mathbf{Z}_l(1-r))})$$

when $L_\mathbf{Q}(\chi, r) = 0$.

**4.8** Some elements in $\text{ann}_{\mathbf{Z}_l[G(\mathbf{Q}(\xi_{ml^{s+1}})^+/\mathbf{Q})]}(\mathcal{H}^1_{\text{cyclo}})$

We shall now describe a method for constructing elements in the annihilator of

$$\mathcal{H}^1_{\text{cyclo}} = \frac{H^1_{\text{ét}}(X_l^+; \mathbf{Z}_l(1-r))}{H^1_{\text{cyclo}}(X_l^+; \mathbf{Z}_l(1-r))}$$

when $0 > r$ is even. The construction will involve the *derivative* of the Dirichlet L-function $L_{\mathbf{Q}}(z, \chi^{-1})$ at $z = r$. In order to present our construction functorially we shall imitate the Stark conjecture [54]. Let

$$Y_r(\mathbf{Q}(\xi_{ml^{s+1}})) = \prod_{\Sigma_{\mathbf{Q}(\xi_{ml^{s+1}})}} (2\pi i)^{-r} \mathbf{Z}$$

with a $G(\mathbf{Q}(\xi_{ml^{s+1}})/\mathbf{Q}) \times G(\mathbf{C}/\mathbf{R})$-action in which $G(\mathbf{C}/\mathbf{R})$ acts on $\Sigma_{\mathbf{Q}(\xi_{ml^{s+1}})}$ and on $(2\pi i)^{-r}$ (the latter action being trivial when $r$ is even). Here, as in §4.2, $\Sigma_{\mathbf{Q}(\xi_{ml^{s+1}})}$ is the set of embeddings of $\mathbf{Q}(\xi_{ml^{s+1}})$ into the complex numbers. The Borel regulator is a Galois equivariant $\mathbf{R}$-linear isomorphism of the form ([5] §3.2)

$$R^r_{\mathbf{Q}(\xi_{ml^{s+1}})^+} : K_{1-2r}(X_l^+) \otimes \mathbf{R} \xrightarrow{\cong} Y_r(\mathbf{Q}(\xi_{ml^{s+1}}))^+ \otimes \mathbf{R} \cong \mathbf{R}[G(\mathbf{Q}(\xi_{ml^{s+1}})^+/\mathbf{Q})]\langle y_r \rangle$$

where, in the notation of ([5] §8), $Y_r(\mathbf{Q}(\xi_{ml^{s+1}}))^+$ is a free $\mathbf{Z}[G(\mathbf{Q}(\xi_{ml^{s+1}})^+/\mathbf{Q})]$-module of rank one with a generator $y_r$.

Define a homomorphism from the complex representation ring to the non-zero complex numbers by sending a character $\chi : G(\mathbf{Q}(\xi_{ml^{s+1}})^+/\mathbf{Q}) \longrightarrow \mathbf{C}^*$ to

$$L^*_{\mathbf{Q}}(r, \chi) = \frac{d}{dz} L_{\mathbf{Q}}(z, \chi)|_{z=r} \cdot \prod_{\substack{p \text{ prime} \\ p | ml, \ (p, f(\chi)) = 1}} (1 - \chi(\sigma_p)^{-1} p^{-r})$$

where $f(\chi)$ denotes the conductor of $\chi$ and $L_{\mathbf{Q}}(z, \chi)$ is the Dirichlet L-function. Since $R^r_{\mathbf{Q}(\xi_{ml^{s+1}})^+}$ is an isomorphism of $\mathbf{R}[G(\mathbf{Q}(\xi_{ml^{s+1}})^+/\mathbf{Q})]$-modules there must exist an isomorphism of $\mathbf{Q}[G(\mathbf{Q}(\xi_{ml^{s+1}})^+/\mathbf{Q})]$-modules of the form

$$f_r : K_{1-2r}(X_l^+) \otimes \mathbf{Q} \xrightarrow{\cong} Y_r(\mathbf{Q}(\xi_{ml^{s+1}}))^+ \otimes \mathbf{Q}.$$

Form the Stark regulator, for each representation $V$ whose contragredient is denoted by $V^\vee$,

$$R(V, f_r) = det((R^r_{\mathbf{Q}(\xi_{ml^{s+1}})^+} \cdot f_r^{-1})^*),$$

where $(R^r_{\mathbf{Q}(\xi_{ml^{s+1}})^+} \cdot f_r^{-1})^*$ is the automorphism of the complex vector space $Hom_{G(\mathbf{Q}(\xi_{ml^{s+1}})^+/\mathbf{Q})}(V^\vee, Y_r(\mathbf{Q}(\xi_{ml^{s+1}}))^+ \otimes \mathbf{C})$ given by composition with $R^r_{\mathbf{Q}(\xi_{ml^{s+1}})^+} \cdot f_r^{-1}$. Following the example of the Stark conjecture [54], we define a complex-valued homomorphism $\mathcal{R}_{f_r}$ given on a finite-dimensional complex representation $V$ by setting

$$\mathcal{R}_{f_r}(V) = \frac{R(V, f_r)}{L^*_{\mathbf{Q}}(r, V)}.$$

If $\overline{\mathbf{Q}}$ is the algebraic closure of the rationals in the complex numbers and $R(G(\mathbf{Q}(\xi_{ml^{s+1}})^+/\mathbf{Q}))$ denotes the representation ring then

$$\mathcal{R}_{f_r} \in Hom_{G(\overline{\mathbf{Q}}/\mathbf{Q})}(R(G(\mathbf{Q}(\xi_{ml^{s+1}})^+/\mathbf{Q})), \overline{\mathbf{Q}}^*).$$

This is seen in the following manner. Since any two choices of $f_r$ differ by multiplication by a unit in $\mathbf{Q}[G(\overline{\mathbf{Q}}/\mathbf{Q}(\xi_{ml^{s+1}})^+/\mathbf{Q})]^*$ the truth of the assertion is independent of $f_r$. Now $K_{1-2r}(X_l^+) \otimes \mathbf{Q}$ is the free $\mathbf{Q}[G(\mathbf{Q}(\xi_{ml^{s+1}})^+/\mathbf{Q})]$-module on the Beilinson element which is denoted by $(-r)!(ml^{s+1})^{-r}b_r(ml^{s+1})$ in the notation of [5]. Therefore we may choose $f_r$ to satisfy

$$f_r((-r)!(ml^{s+1})^{-r}b_r(ml^{s+1})) = y_r.$$

Let $G$ be a finite abelian group. Then there exists a homomorphism

$$\lambda_G : Hom(R(G), \mathbf{C}^*) \xrightarrow{\cong} \mathbf{C}[G]^*$$

given by

$$\lambda_G(h) = \sum_{\chi \in \hat{G}} h(\chi)e_\chi$$

where $e_\chi = |G|^{-1} \sum_{g \in G} \chi(g)g^{-1} \in \mathbf{Q}(\chi)[G]$. Furthermore $\lambda_G$ restricts to give an isomorphism of the form

$$\lambda_G : Hom_{\Omega_{\mathbf{Q}}}(R(G), \overline{\mathbf{Q}}^*) \xrightarrow{\cong} \mathbf{Q}[G]^*.$$

Now we evaluate $\lambda_{G(\mathbf{Q}(\xi_{ml^{s+1}})^+/\mathbf{Q})}(\mathcal{R}_{f_r})$. Clearly

$$\lambda_{G(\mathbf{Q}(\xi_{ml^{s+1}})^+/\mathbf{Q})}(V \mapsto L_{\mathbf{Q}}^*(r, V)) = \sum_\chi L_{\mathbf{Q}}^*(r, \chi)e_\chi.$$

On the other hand the subspace of $\mathbf{C}[G(\mathbf{Q}(\xi_{ml^{s+1}})^+/\mathbf{Q})]\langle y_r \rangle$ on which $G(\mathbf{Q}(\xi_{ml^{s+1}})^+/\mathbf{Q})$ acts via $\chi^{-1}$ is $\mathbf{C}e_{\chi^{-1}} \cdot y_r$. Hence

$$\lambda_{G(\mathbf{Q}(\xi_{ml^{s+1}})^+/\mathbf{Q})}(V \mapsto R(V, f_r)) = \sum_\chi \mu_{\chi^{-1}} e_\chi$$

$R^r_{\mathbf{Q}(\xi_{ml^{s+1}})^+} \cdot f_r^{-1}$ is given by multiplication by $\mu_\chi$ on $\mathbf{C}e_\chi \cdot y_r$. However, as recapitulated in ([5] §4.2), Beilinson proved ([37] Part I Theorem 4.3(ii) and Part II Theorem 1.1) that

$$R^r_{\mathbf{Q}(\xi_{ml^{s+1}})^+}(f_r^{-1}(y_r)) = R^r_{\mathbf{Q}(\xi_{ml^{s+1}})^+}((-r)!(ml^{s+1})^{-r}b_r(ml^{s+1}))$$

$$= \sum_\chi L_{\mathbf{Q}}^*(r, \chi^{-1})e_\chi y_r$$

so that $\mu_{\chi^{-1}} = L_{\mathbf{Q}}^*(r, (\chi^{-1})^{-1}) = L_{\mathbf{Q}}^*(r, \chi)$ and therefore $\lambda_{G(\mathbf{Q}(\xi_{ml^{s+1}})^+/\mathbf{Q})}(\mathcal{R}_{f_r}) = 1 \in \mathbf{Q}[G]^*$, as required.

Now choose any $f_r$ and define

$$\mathcal{I}_{f_r} = \{\alpha \in \mathbf{Q}[G(\mathbf{Q}(\xi_{ml^{s+1}})^+/\mathbf{Q})] \mid \alpha \cdot f_r(K_{1-2r}(X_l^+)) \subseteq Y_r(\mathbf{Q}(\xi_{ml^{s+1}}))^+\}.$$

Hence $\mathcal{I}_{f_r}$ is a $\mathbf{Z}[G(\mathbf{Q}(\xi_{ml^{s+1}})^+/\mathbf{Q})]$-submodule of $\mathbf{Q}[G(\mathbf{Q}(\xi_{ml^{s+1}})^+/\mathbf{Q})]$.

**Theorem 4.9**

In the notation of §3.7 the fractional ideal of $\mathbf{Q}[G(\mathbf{Q}(\xi_{ml^{s+1}})^+/\mathbf{Q})]$

$$\mathcal{J}_r = \mathcal{I}_{f_r} \cdot \lambda_{G(\mathbf{Q}(\xi_{ml^{s+1}})^+/\mathbf{Q})}(\mathcal{R}_{f_r})^{-1}$$

is independent of the choice of $f_r$ and for each odd prime $l$

$$\mathcal{J}_r \bigcap \mathbf{Z}_l[G(\mathbf{Q}(\xi_{ml^{s+1}})^+/\mathbf{Q})] \subseteq \mathrm{ann}_{\mathbf{Z}_l[G(\mathbf{Q}(\xi_{ml^{s+1}})^+/\mathbf{Q})]}(H^2_{\mathrm{\acute{e}t}}(X_l^+; \mathbf{Z}_l(1-r)))$$

when $0 > r$ is even.

**Proof**

Any other choice for $f_r$ has the form $uf_r$ for some unit $u$ in the rational group-ring. This changes $\mathcal{I}_{f_r}$ to $u^{-1}\mathcal{I}_{f_r}$ and changes $\lambda_{G(\mathbf{Q}(\xi_{ml^{s+1}})^+/\mathbf{Q})}(\mathcal{R}_{f_r})$ to $u^{-1}\lambda_{G(\mathbf{Q}(\xi_{ml^{s+1}})^+/\mathbf{Q})}(\mathcal{R}_{f_r})$ so that

$$\mathcal{I}_{uf_r}\lambda_{G(\mathbf{Q}(\xi_{ml^{s+1}})^+/\mathbf{Q})}(\mathcal{R}_{uf_r})^{-1} = u^{-1}\mathcal{I}_{f_r}u\lambda_{G(\mathbf{Q}(\xi_{ml^{s+1}})^+/\mathbf{Q})}(\mathcal{R}_{f_r})^{-1}$$

$$= \mathcal{I}_{f_r}\lambda_{G(\mathbf{Q}(\xi_{ml^{s+1}})^+/\mathbf{Q})}(\mathcal{R}_{f_r})^{-1},$$

as required.

For the second part, choose the canonical $f_r$ as in §4.8. Let $x \in H^1_{\mathrm{\acute{e}t}}(X_l^+; \mathbf{Z}_l(1-r))$ and $\alpha \in \mathcal{J}_r \bigcap \mathbf{Z}_l[G(\mathbf{Q}(\xi_{ml^{s+1}})^+/\mathbf{Q})]$. By Theorem 4.6 we must show that $\alpha x \in H^1_{\mathrm{cyclo}}(X_l^+; \mathbf{Z}_l(1-r))$. The image of the Chern class

$$c_{1-r,1}: K_{1-2r}(X_l^+) \longrightarrow H^1_{\mathrm{\acute{e}t}}(X_l^+; \mathbf{Z}_l(1-r))$$

(denoted by $c^r_{\mathbf{Q}(\xi_{ml^{s+1}})^+}$ in ([5] Lemma 8.16)) is dense. Choose $x_m \in K_{1-2r}(X_l^+)$ so that $\lim_{\overrightarrow{m}} c_{1-r,1}(x_m) = x$. Then $f_r(\alpha x_m) \in Y_r(\mathbf{Q}(\xi_{ml^{s+1}}))^+$ which equals $f_r(\mathbf{Z}[G(\mathbf{Q}(\xi_{ml^{s+1}})^+/\mathbf{Q})]\langle(-r)!(ml^{s+1})^{-r}b_r(ml^{s+1})\rangle)$ so that

$$\alpha x_m \in \mathbf{Z}[G(\mathbf{Q}(\xi_{ml^{s+1}})^+/\mathbf{Q})]\langle(-r)!(ml^{s+1})^{-r}b_r(ml^{s+1})\rangle.$$

However ([5] Lemma 8.16; [22] Theorem 6.4 (proof))) $c_{1-r,1}((-r)!(ml^{s+1})^{-r}b_r(ml^{s+1}))$ is equal to the element $\eta_m(r)$ of §4.2 and §4.5, which generates the cyclotomic elements. Hence $\lim_{\overrightarrow{m}} \alpha c_{1-r,1}(x_m) = \alpha x \in H^1_{\mathrm{cyclo}}(X_l^+; \mathbf{Z}_l(1-r))$, as required.
□

# 5 Other invariants in $K_0(\mathbf{Z}_l[G], \mathbf{Q}_l)$

**5.1** *Relatively abelian number field extensions*

Let $l$ be an odd prime and let $L/K$ be a finite Galois extension of number fields with abelian Galois group $G(L/K)$. Assume that $L$ is totally real. Let $K_\infty/K$ denote the cyclotomic $\mathbf{Z}_l$-extension of $K$ and suppose that $L \cap K_\infty =$ $K$. Set $L_\infty = K_\infty L$ and let $M_\infty$ denote the maximal abelian pro-$l$-extension of $L_\infty$ in which only primes above $l$ ramify. Therefore we have an extension of Galois groups of the form

$$X_\infty = G(M_\infty/L_\infty) \longrightarrow G(M_\infty/K) \longrightarrow G(L_\infty/K) \cong G(L/K) \times G(K_\infty/K)$$

and $G(K_\infty/K) \cong \mathbf{Z}_l$. It is usual to write $\Gamma$ for $G(K_\infty/K)$. Hence $X_\infty$ is a $\mathbf{Z}_l[[G(L_\infty/K)]]$-module and, in particular, it is a $\Gamma$-module. Let $\overline{\mathbf{Q}}_l$ denote an algebraic closure of the $l$-adics. For each character $\eta : G(L/K) \longrightarrow \overline{\mathbf{Q}}_l^*$ of order prime to $l$ the $\eta$-eigenspace

$$X_\infty^\eta = \{x \in X_\infty \otimes \mathbf{Z}_l[\eta] \mid g(x) = \eta(g)x \ for\ all\ g \in G(L/K)\}$$

is also a $\Gamma$-module. Let $\mu_\eta$ denote the Iwasawa $\mu$-invariant of $X^\eta$ ([38] p.604 and p.655 et seq, [50], [59]).

Suppose that $\mu_\eta = 0$ for all $l$-adic characters $\eta : G(L/K) \longrightarrow \overline{\mathbf{Q}}_l^*$ of order prime to $l$. This condition holds when $K$ is the field of rational numbers by [17]. In [42] a perfect complex of $\mathbf{Z}_l[G(L_\infty/K)]$-modules is constructed which generalises that constructed in [5] for the cyclotomic case. Unfortunately, unless one assumes the Leopoldt Conjecture ([58] p.71) it does not seem possible to establish a descent lemma similar to that of ([5] §8.1). On the other hand, if we make this assumption then Theorem 1.6 should extend to this situation.

**5.2** *Vanishing cycles*

In [23] the following situation is considered: $\mathcal{O}_k$ is an excellent, henselian discrete valuation ring with fraction field, $k$, and algebraically closed residue field, $F$. Then $\mathcal{O}_k[T]$ is a local ring with maximal ideal $\mathcal{M} =< T, P = ker(\mathcal{O}_k \longrightarrow F) >$. Let $A = \mathcal{O}_k\{T\}$ denote the henselianisation of $\mathcal{O}_k[T]$ at $\mathcal{M}$ ([36] p.36). We also write $\mathcal{P} \triangleleft A$ for the kernel of $\mathcal{O}_k\{T\} \longrightarrow F\{T\}$ so that $\mathcal{P} \cap \mathcal{O}_k = P$. If $X = \mathrm{Spec}(A) - \{\mathcal{M}, \mathcal{P}\}$ and $P_1, P_2, \ldots, P_t$ is a non-empty finite set of primes of $X$ different from $\mathcal{M}$ and $\mathcal{P}$ then we have

$$U = \mathrm{Spec} A - \{\mathcal{M}, \mathcal{P}, P_1, P_2, \ldots, P_t\} \stackrel{u}{\hookrightarrow} X \stackrel{i}{\hookleftarrow} Y = \{P_1, P_2, \ldots, P_t\}$$

where $u$ is an open immersion and $i$ is a closed immersion.

Let $\Lambda = \mathbf{F}_q$ be a finite field of characteristic $l$ different from that of $F$ and let $\mathcal{F}$ be a locally constant sheaf of $\Lambda$-modules of finite rank on $U_{\acute{e}t}$. Fix $\overline{k}$, an algebraic closure of $k$. Then ([23] Theorem 6.7) gives a formula, generalising that of Deligne [28], for the $\Lambda$-dimension of the space of vanishing cycles

$$H^q(X \times_{\mathrm{Spec}(k)} \mathrm{Spec}(\overline{k}); u_!(\mathcal{F}))$$

in terms of the Swan conductor. In general this group is zero except when $q = 0, 1$ ([21] p.14) and is also zero when $q = 0$ in the above situation because $U \neq X$ ([23] p.631).

There is an action on this cohomology of the absolute Galois group $G(\overline{k}/k)$ on $H^q(X \times_{\text{Spec}(k)} \text{Spec}(\overline{k}); u_!(\mathcal{F}))$ and the fundamental question posed in [23] is to determine this Galois action. Particularly important is the case when $q = 1$ and $U \neq X$. The following construction is a small step towards studying this problem.

Suppose that $k$ is a complete discrete valuation field of finite cohomological dimension. For example, $k$ might be an $n$-dimensional local field (see [18], [19]). Assume also that $U \neq X$. For each finite Galois extension, $k \subseteq L \subset \overline{k}$ there exists a perfect complex [51], whose $s$-th cohomology group is equal to $H^{s+1}(X \times_{\text{Spec}(k)} \text{Spec}(L); u_!(\mathcal{F}))$, which defines an Euler characteristic in the relative K-group

$$\chi(X \times_{\text{Spec}(k)} \text{Spec}(L); u_!(\mathcal{F})) = [P^{od}, \Gamma, P^{ev}] \in K_0(\mathbf{Z}_l[G(L/k)], \mathbf{Q}_l).$$

When $G(L/k)$ is abelian the class $[P^{od}, \Gamma, P^{ev}]$ corresponds to

$$det(\Gamma) \in \frac{\mathbf{Q}_l[G(L/k)]^*}{\mathbf{Z}_l[G(L/k)]^*} \cong K_0(\mathbf{Z}_l[G(L/k)], \mathbf{Q}_l).$$

The formula of ([23] Theorem 6.7) suggests that there should an expression for $det(\Gamma)$ in terms of Swan conductors. If, in addition, $k$ is a one-dimensional local field then $H^s(X \times_{\text{Spec}(k)} \text{Spec}(L); u_!(\mathcal{F}))$ can only be possibly non-zero for $s = 1, 2, 3$. There are generalisations of Theorem 2.4 to the case of three non-zero homology groups, but the ones which I have currently do not give very efficient annihilator/Fitting ideal relations. For example, the expected form of one of these relations would predict that some power of the annihilators of $H^1(X \times_{\text{Spec}(k)} \text{Spec}(L); u_!(\mathcal{F}))$ and $H^3(X \times_{\text{Spec}(k)} \text{Spec}(L); u_!(\mathcal{F}))$ times an expression involving Swan conductors lies inside the annihilator ideal of $H^2(X \times_{\text{Spec}(k)} \text{Spec}(L); u_!(\mathcal{F}))$.

### 5.3 *Geometrical application to surfaces*

Let $X$ be a smooth projective, geometrically connected surface over a finite field of characteristic $p$ and assume that $Br(X)$ is finite (as is widely believed to be true). Let $G$ be a finite group acting on $X$ in such a way that $\pi : X \longrightarrow X/G$ is étale. The assumption that $\pi$ is étale ensures that ([9] Proposition 3.2) that there exists a bounded cochain complex of cohomologically trivial $\mathbf{Z}[G]$-modules $C^*$ such that $H^i(X; \mathbf{G}_m)$ and $H^i(C^*)$ are isomorphic $\mathbf{Z}[G]$-modules for each $i$. In [9], using the results of ([29], [30], [31], [32]), an Euler characteristic is constructed in $K_0(\mathbf{Z}[G])$ whose image in the Grothendieck group of finite $\mathbf{Z}[G]$-modules is given in terms of the finite parts of $H^*(X; \mathbf{G}_m)$ and the intersection pairing on $H^1(X; \mathbf{G}_m)$. In [6] this Euler characteristic is lifted canonically to the Euler characteristic $[\tilde{P}^{od}, Z, \tilde{P}^{ev}] \in K_0(\mathbf{Z}[G], \mathbf{Q})$ and associated to a perfect complex $\tilde{P}^*$. If $G$ is abelian the first obstacle to applying Theorem 2.4 to $\tilde{P}^*$ is the fact that the cohomology groups are not finite. This may be remedied by choosing chain map of a perfect complexes $f : P_1^* \to \tilde{P}^*$, which induces isomorphisms on $H^* \otimes \mathbf{Q}_l$, and replacing $\tilde{P}^*$ by the mapping cone of $f$. Unfortunately the mapping cone will usually have more than two non-zero (finite) cohomology groups so that, as in §5.2, we require an efficient version of Theorem 2.4 applicable to a bounded perfect complex with finite homology groups.

On the positive side when $G$ is abelian the computations of ([9] Theorem 3.9) suggests that, for each prime $l \neq p$, $det(Z) \in \mathbf{Q}_l[G]^*/\mathbf{Z}_l[G]^*$ is given in terms of the Artin L-function of $X \to X/G$.